\documentclass{article}
\usepackage{amsthm,amsmath,amssymb,latexsym,amscd}

\renewcommand{\a}{\alpha}
\renewcommand{\b}{\beta}
\newcommand{\g}{\frak{g}}
\newcommand{\h}{\frak{h}}

\newcommand{\A}{\mathcal{A}}

\newcommand{\T}{\mathcal{T}}
\newcommand{\I}{\mathcal{I}}

\newcommand{\Hom}{\mathrm{Hom}}

\newcommand{\p}{\prime}
\newcommand{\mb}{\mathbf}
\renewcommand{\c}{\circ}
\newcommand{\ot}{\otimes}

\newcommand{\what}{\widehat}

\newtheorem{definition}{Definition}[section]
\newtheorem{lemma}[definition]{Lemma}
\newtheorem{proposition}[definition]{Proposition}
\newtheorem{theorem}[definition]{Theorem}
\newtheorem{corollary}[definition]{Corollary}

\newtheorem{example}[definition]{Example}
\date{}

\begin{document}

\title{Twisting on associative algebras\\
and\\
Rota-Baxter type operators}

\author{Kyousuke UCHINO}

\maketitle

\noindent
{\footnotesize
Supported by Keio University COE program}.\\
{\footnotesize
Science University of Tokyo,
Wakamiya 26 Shinjyuku Tokyo Japan}.\\
{\footnotesize \textbf{e-mail}:
K\underline{ }Uchino[at]oct.rikadai.jp}.
\medskip\\
\noindent
{\footnotesize \textbf{MSC (2000)}:
13D03, 16S80, 17B62.}\\
\noindent
{\footnotesize \textbf{Keywords}: deformation theory, twisting,
Rota-Baxter operators, Reynolds operators, Nijenhuis operators}.

\abstract
{
We will introduce an operation ``twisting" on Hochschild complex
by analogy with Drinfeld's twisting operations.
By using the twisting and derived bracket construction,
we will study differential graded Lie
algebra structures associated with bi-graded Hochschild complex.
We will show that Rota-Baxter type operators
are solutions of Maurer-Cartan equations.
As an application of twisting, we will give a construction
of associative Nijenhuis operators.
}

\section{Introduction.}

In \cite{D}, Drinfeld introduced an operation ``twisting",
motivated by the study of quasi-Lie bialgebras and
quasi-Hopf algebras.
The twisting operations provide a method
of analyzing Manin triples.
In the context of Poisson geometry,
they gave the detailed study of twisting operations
(see Kosmann-Schwarzbach \cite{Kos1,Kos3}
and Roytenberg \cite{Roy1,Roy2}).
We shortly describe twisting operations.
We consider a graded commutative algebra,
$\bigwedge^{\cdot}(V\oplus V^{*})$,
where $V$ is a vector space over $\mathbb{R}$,
$V^{*}$ is the dual space of $V$.
The graded algebra has a graded
Poisson bracket defined by
$\{V,V\}=\{V^{*},V^{*}\}:=0$ and
$\{V,V^{*}\}:=\langle{V,V^{*}}\rangle$.
By definition, a {\em structure}
in the graded Poisson algebra is an element $\Theta$ in
$\bigwedge^{3}(V\oplus V^{*})$ satisfying a
Maurer-Cartan equation $\{\Theta,\Theta\}=0$.
It is known that
the structure $\Theta$ is an invariant
Lie algebra structure on $V\oplus V^{*}$.
The structures are closely related with
(quasi-)Lie bialgebra structures.
A Lie bialgebra structure is defined as
a pair of tensors $(\nu_{1},\nu_{2})$ such that
$\Theta_{12}:=\nu_{1}+\nu_{2}$ is a structure
in above sense,
where $\nu_{1}\in(\bigwedge^{2}V^{*})\ot V$
and $\nu_{2}\in V^{*}\ot\bigwedge^{2}V$.
When $(\nu_{1},\nu_{2})$ is a structure of Lie bialgebra,
the total space $(V\oplus V^{*},\Theta_{12})$
is called a Drinfeld double.
Let $r$ be an element in $V\wedge V$.
By definition,
the twisting of a structure $\Theta$
by $r$ is a canonical transformation;
$$
\Theta^{r}:=exp(X_{r})(\Theta),
$$
where $X_{r}$ is a Hamiltonian vector field
$X_{r}:=\{-,r\}$ and
$\Theta^{r}$ is the result of twisting.
Several interesting information is riding on
the orbits of twisting operations.
We recall a basic proposition.
Let $(\nu_{1},0)$ be a structure of
Lie bialgebra such that $\nu_{2}=0$.
Then the Drinfeld double is the space $V\oplus V^{*}$
with the structure $\Theta_{1}:=\nu_{1}$.
If $r$ is a solution of a Maurer-Cartan equation
(or classical Yang-Baxter equation)
$$
[r,r]=0,
$$
then a pair $(\nu_{1},\{\nu_{1},r\})$ is a Lie bialgebra
structure and the {\em double}, $\nu_{1}+\{\nu_{1},r\}$,
is equal with the result of twisting $\Theta^{r}_{1}$,
where $[r,r]:=\{\{\nu_{1},r\},r\}$.
Conversely, the Maurer-Cartan condition of $r$
is characterized by this proposition.
\medskip\\
\indent
The aim of this note is to construct the theory of twisting
on associative algebras along the philosophy
and construction in \cite{Kos3} and \cite{Roy1}.
At first, we will define a twisting operation
in the category of associative algebras.
The twisting operation is defined by
using only a canonical bigraded system
of the graded Poisson algebra
$\bigwedge^{\cdot}(V\oplus V^{*})$.
Hence, given a suitable bigraded Lie algebra,
one can define a twisting like operation
on the bigraded Lie algebra.
We consider a Hochschild complex
$C^{*}(\T):=\Hom(\T^{\ot *},\T)$, where $\T$ is a vector space
decomposed into two subspaces $\T:=\A_{1}\oplus\A_{2}$.
In Section 2, we will introduce
a canonical bigraded Lie algebra system on
$C^{*}(\A_{1}\oplus\A_{2})$.
The graded Lie bracket is given by Gerstenhaber's
bracket product.
Our {\em structures}, $\theta$,
are defined as associative structures
on $\A_{1}\oplus\A_{2}$, i.e.,
$\theta$ is a 2-cochain in $C^{2}(\A_{1}\oplus\A_{2})$
and $t_{1}*t_{2}:=\theta(t_{1}\ot t_{2})$
is associative for any $t_{1},t_{2}\in\A_{1}\oplus\A_{2}$.
For a given 1-cochain $H:\A_{2}\to\A_{1}$,
we define a twisting operation
by the same manner with classical one,
$$
\theta^{H}:=exp(X_{\what{H}})(\theta),
$$
where $\what{H}$ is the image of the natural map
$C^{*}(\A_{2},\A_{1})\hookrightarrow C^{*}(\A_{1}\oplus\A_{2})$
and $X_{\what{H}}$ is an analogy of Hamiltonian
vector field defined by $X_{\what{H}}:=\{-,\what{H}\}$,
where $C^{*}(\A_{2},\A_{1}):=\Hom(\A_{2}^{\ot *},\A_{1})$.
We will see that $\theta$ is decomposed into the unique
4 substructures,
$$
\theta=\hat{\phi}_{1}+
\hat{\mu}_{1}+\hat{\mu}_{2}+\hat{\phi}_{2}.
$$
The twisting operation is completely determined
by transformation rules of the 4 substructures.
In Section 4, we will give explicit formulas
of the transformation rules
(Theorem \ref{maintheorem}).\\
\indent
We consider the case of
$\hat{\phi}_{1}=\hat{\phi}_{2}=0$.
In this case, $\A_{1}$ and $\A_{2}$ are both subalgebras of
the associative algebra $(\A_{1}\oplus\A_{2},\theta)$.
Such a triple $(\A_{1}\oplus\A_{2},\A_{1},\A_{2})$ is called
an associative twilled algebra, simply, twilled algebra
(Carinena and coauthors \cite{CGM}).
When a Lie algebra is decomposed into two subalgebras,
it is called a {\em twilled} Lie algebra (\cite{Kos2}),
or called a twilled extension in \cite{KM},
or a {\em double} Lie algebra in \cite{LW}.
This concept is used in order to construct
integrable Hamiltonian systems
(Adler-Kostant-Symes theorem).
The notion of associative twilled algebra
is considered as an associative version
of the classical one.
In \cite{CGM},
they studied associative twilled algebras
from the point of view of quantization.
In Section 3, we will give the detailed
study for twilled algebras.
By derived bracket construction in \cite{Kos2},
a twilled algebra structure on $\A_{1}\oplus\A_{2}$
induces a differential graded Lie algebra
(shortly, dg-Lie algebra) structure on $C^{*}(\A_{2},\A_{1})$
(see Proposition \ref{dglie}).
So we can consider a deformation theory
on the induced dg-Lie algebra.
We consider a Maurer-Cartan equation in the dg-Lie algebra,
$$
dR+\frac{1}{2}[R,R]=0.
$$
We can find a solution $R$ in Rota-Baxter algebra theory.
Let $(\A,R)$ be an arbitrary associative algebra equipped
with an operator $R:\A\to\A$.
The operator $R$ is called a Rota-Baxter operator,
if $R$ satisfies an identity (so-called Rota-Baxter identity),
\begin{equation*}
R(x)R(y)=R(R(x)y+xR(y))+qR(xy),
\end{equation*}
where $q\in\mathbb{K}$ is a scalar (called a weight).
Rota-Baxter operators have been studied in combinatorics
(see Rota \cite{Rot1,Rot2}).
In this note we do not study the combinatorial problem,
because it is beyond our aim.
Now $\A\oplus\A$ has a natural twilled algebra structure,
and then $C^{*}(\A,\A)$ has a dg-Lie algebra structure.
In Section 5.1, we will show that
$R$ is a Rota-Baxter operator if and only if
$R$ is a solution of the Maurer-Cartan equation.\\
\indent
In Section 6, we will give an application of our construction.
We recall the notion of associative Nijenhuis operator (\cite{CGM}).
Let $N:\A\to\A$ be a linear map on an associative algebra $\A$.
The operator $N$ is called an
{\em associative Nijenhuis operator},
if it satisfies an associative version of
classical Nijenhuis condition,
$$
N(x)N(y)=N(N(x)y+xN(y))-N^{2}(xy),
$$
where $x,y\in\A$.
They showed that
a deformed multiplication,
$x\times_{N}y:=N(x)y+xN(y)-N(xy)$,
is a new associative multiplication on $\A$
and it is compatible with original one.
In this sense, an associative Nijenhuis operator
induces a quantum bihamiltonian system (see \cite{CGM}).
We will give a construction of associative
Nijenhuis operators by analogy
with Poisson-Nijenhuis geometry.\\
\indent
We recall a theorem of Vaisman \cite{V}.
Let $(V,P)$ be a Poisson manifold equipped with
a Poisson structure tensor $P$, i.e., $P$ is a solution of
a Maurer-Cartan equation,
$$
\frac{1}{2}[P,P]=0,
$$
where the bracket product is a graded Lie bracket
(called Schouten-Nijenhuis bracket).
Since the Poisson structure is a $(2,0)$-tensor,
it is identified with a bundle map $P:T^{*}V\to TV$.
The Poisson bundle map induces a Lie algebroid
structure on the cotangent bundle $T^{*}V$,
i.e., the space of sections of $\bigwedge^{\cdot}T^{*}V$
has a certain graded Lie bracket $\{,\}_{P}$.
He showed that if a 2-form $\omega$ is a solution of
the strong Maurer-Cartan equation,
$d\omega=\{\omega,\omega\}_{P}=0$, then the bundle map
$N:=P\omega:TV\to TV$ is a Nijenhuis tensor and the pair $(P,N)$ is
a compatible pair, or a Poisson-Nijenhuis
structure in the sense of \cite{Kos5}.
This compatibility implies that the bundle map
$NP:T^{*}V\to TV$ is
a Poisson structure bundle map and $P+t NP$ is
a one parameter family of Poisson structures.\\
\indent
We will show a similar theorem to Vaisman's theorem.
First of all, we need Rota-Baxter {\em type} operators
as substitutes for Poisson structures.
Let $\A$ be an associative algebra and $M$ an $\A$-bimodule,
and let $\pi:M\to\A$ be a linear map.
The linear map $\pi$ is called
a {\em generalized Rota-Baxter operator}
of weight 0, or shortly GRB (\cite{U1}),
if $\pi$ is a solution of
$$
\pi(m)\pi(n)=\pi(\pi(m)\cdot n+m\cdot \pi(n)),
\eqno{(GRB)}
$$
where $m,n\in M$ and $\cdot$ is the bimodule action.
When $M=\A$ as a canonical bimodule, (GRB)
reduces to a classical Rota-Baxter identity of weight zero.
We consider a semidirect product algebra
$(\T:=A\ltimes M,\hat{\mu})$,
where $\hat{\mu}$ is the associative multiplication
of $\A\ltimes M$.
The Hochschild complex $C^{*}(\A\ltimes M)$
becomes a dg-Lie algebra by Gerstenhaber bracket
and the coboundary map
$d_{\hat{\mu}}:=\{\hat{\mu},-\}$.
We define, due to \cite{Kos2},
a second bracket product on $C^{*}(\A\ltimes M)$ by
$$
[f,g]_{\hat{\mu}}:=(-1)^{|f|-1}\{\{\hat{\mu},f\},g\}.
$$
Here the new bracket is a graded Lie bracket
on $C^{*}(M,\A)\subset C^{*}(\A\ltimes M)$.
One can show that $\pi$ is a generalized Rota-Baxter
operator if and only if it is a solution of
the Maurer-Cartan equation
$$
\frac{1}{2}[\hat{\pi},\hat{\pi}]_{\hat{\mu}}=0,
$$
where $\hat{\pi}$ is the image of the natural map
$C^{1}(M,\A)\hookrightarrow C^{1}(\A\ltimes M)$,
$\pi\mapsto\hat{\pi}$.\\
\indent
Now, given a generalized Rota-Baxter operator
$\pi:M\to\A$, $M$ becomes an associative algebra,
where the associative multiplication on $M$ is
given by a structure $\{\hat{\mu},\hat{\pi}\}$.
The associativity of $\{\hat{\mu},\hat{\pi}\}$
is followed from $[\hat{\pi},\hat{\pi}]_{\hat{\mu}}=0$.
We denote the associative algebra by $M_{\pi}$.
One can show that $M_{\pi}\oplus\A$ has a twilled algebra
structure. Thus a dg-Lie algebra structure,
$(d_{\hat{\mu}},[,]_{\{\hat{\mu},\hat{\pi}\}})$,
is induced on $C^{*}(\A,M_{\pi})$.
By analogy with Vaisman's theorem, we assume that
$\Omega:\A\to M$ is a solution of the strong
Maurer-Cartan equation in $C^{*}(\A,M_{\pi})$,
$$
d_{\hat{\mu}}\hat{\Omega}=
[\hat{\Omega},\hat{\Omega}]_{\{\hat{\mu},\hat{\pi}\}}=0,
$$
where $d_{\hat{\mu}}$ is the Hochschild coboundary on $C^{*}(\A,M)$
and $\hat{\Omega}$ is defined by the similar manner with $\hat{\pi}$.
Then we can show that a linear endomorphism $N:=\pi\Omega:\A\to\A$
is an associative Nijenhuis operator and the pair $(\pi,N=\pi\Omega)$
is compatible (see Proposition \ref{constN}).
This proposition can be considered as an associative version of
Vaisman's result.
\medskip\\
\noindent
\textbf{Acknowledgements.}
The author wishes to thank very much the referees.
He is greatly indebted to them
for their numerous suggestion.
Finally, he would like to thank very
much Professors Jean-Louis Loday,
Yoshiaki Maeda and Akira Yoshioka
for helpful comments and encouragement.

\section{Cochain calculus.}

In this section, we will define a bigraded Lie algebra
structure on Hochschild complex $C^{*}(\A_{1}\oplus\A_{2})$.
In the following, we assume that
the characteristic of a ground field
$\mathbb{K}$ is zero and that
$\mathbb{Q}$ is included in $\mathbb{K}$.

\subsection{Gerstenhaber brackets.}

We recall Gerstenhaber's bracket product.
Let $V$ be a vector space over $\mathbb{K}$.
We consider the space of cochains
$\g(V):=\bigoplus_{n\in\mathbb{N}}C^{n}(V)$,
where $C^{n}(V)=C^{n}(V,V):=\Hom_{\mathbb{K}}(V^{\ot n},V)$.
By definition, the degree of $f\in\g(V)$ is $|f|$,
if $f$ is in $C^{|f|}(V)$.
For any $f\in C^{|f|}(V)$ and
$g\in C^{|g|}(V)$, we define a product,
$$
f\bar{\c}g:=\sum^{|f|}_{i=1}(-1)^{(i-1)(|g|-1)}f\c_{i}g,
$$
where $\c_{i}$ is the composition of maps defined by
$$
f\c_{i}g(b_{1},...,b_{|f|+|g|-1})=
f(b_{1},...,b_{i-1},g(b_{i},...,b_{i+|g|-1}),b_{i+|g|}...,b_{|f|+|g|-1}).
$$
The degree of $f\bar{\c}g$ is $|f|+|g|-1$.
The Gerstenhaber bracket, or shortly, G-bracket on $\g(V)$
is defined as a graded commutator,
$$
\{f,g\}:=f\bar{\c}g-(-1)^{(|f|-1)(|g|-1)}g\bar{\c}f.
$$
We recall two fundamental identities:\\
(I) graded commutativity,
$$
\{f,g\}=-(-1)^{(|f|-1)(|g|-1)}\{g,f\},
$$
(II) graded Jacobi identity,
\begin{multline*}
(-1)^{(|f|-1)(|h|-1)}\{\{f,g\},h\}+
(-1)^{(|h|-1)(|g|-1)}\{\{h,f\},g\}+\\
(-1)^{(|g|-1)(|f|-1)}\{\{g,h\},f\}=0,
\end{multline*}
where $h\in C^{|h|}(V)$.
The above graded Jacobi identity is equivalent with
$$
\{f,\{g,h\}\}=\{\{f,g\},h\}+(-1)^{(|f|-1)(|g|-1)}\{g,\{f,h\}\}.
\eqno{(II^{\p})}
$$
(II$^{\p}$) is called a graded Leibniz
identity, or sometimes called, a graded Loday identity.
\medskip\\
\noindent
\textbf{Graded Lie algebras.}
Let $\g$ be a graded vector space equipped with
a binary multiplication $\{,\}$ of degree $0$.
When the bracket product satisfies the following
two conditions (\ref{defdgl1}) and (\ref{defdgl2}),
$\g$ is called a graded Lie algebra.
\begin{eqnarray}
\label{defdgl1}
\{f,g\}&=&-(-1)^{deg(f)deg(g)}\{g,f\},\\
\label{defdgl2}
\{f,\{g,h\}\}&=&\{\{f,g\},h\}+(-1)^{deg(f)deg(g)}\{g,\{f,h\}\},
\end{eqnarray}
where $f,g,h,\in\g$ and $deg(-)$ is the degree of $\g$.
The cochain complex $\g(V)$ is a graded Lie algebra of $deg(f):=|f|-1$.
A graded Lie algebra $\g$ is called a differential graded Lie algebra
(dg-Lie algebra), if $\g$ has a square zero derivation $d$
of degree $+1$ satisfying,
\begin{equation}\label{defdgl3}
d\{f,g\}=\{d{f},g\}+(-1)^{deg(f)}\{f,d{g}\}.
\end{equation}
\noindent
\textbf{Associative structures.}
It is well-known that $S\in C^{2}(V)$ is an associative structure
if and only if it is a solution of Maurer-Cartan equation,
$\{S,S\}=0$.
If $S$ is an associative structure, then
$d_{S}(f):=\{S,f\}$ is a coboundary map of
Hochschild complex $(C^{*}(V),d_{S})$, and then
$(\g(V),d_{S})$ becomes a dg-Lie algebra.
\medskip\\
\noindent\textbf{Derived brackets.}
Let $\g$ be a dg-Lie algebra.
We define a new bracket product by
$$
[f,g]_{d}:=(-1)^{deg(f)}\{df,g\}.
$$
The new bracket is called
a {\em derived bracket} (\cite{Kos2}).
It is well-known that the derived bracket
is a graded Leibniz bracket, i.e.,
(\ref{defdgl2}) holds up to degree shift.
Remark that the derived bracket is not
graded commutative in general.
We recall a basic lemma.
\begin{lemma}\label{balma}
(\cite{Kos2})
Let $\g$ be a dg-Lie algebra, and let $\h\subset\g$
be an abelian subalgebra, i.e., $\{\h,\h\}=0$.
We define a new degree (derived degree)
by $deg_{d}(h):=deg(h)+1$.
If the derived bracket is closed on $\h$,
then $(\h,deg_{d},[,]_{d})$ is a graded Lie algebra.
\end{lemma}

\subsection{Lift and Bidegree.}

Let $\A_{1}$ and $\A_{2}$ be vector spaces,
and let $c:\A^{\ot n}_{2}\to\A_{1}$ be a linear map,
or a cochain in $C^{n}(\A_{2},\A_{1})$.
We can construct a cochain
$\hat{c}\in C^{n}(\A_{1}\oplus\A_{2})$ by
$$
\hat{c}\Big((a_{1},x_{1})\ot...\ot(a_{n},x_{n})\Big)
:=(c(x_{1},...,x_{n}),0).
$$
In general, for a given multilinear map
$f:\A_{i(1)}\ot\A_{i(2)}\ot...\ot\A_{i(n)}\to\A_{j}$,
$i(1),...,i(n),j\in\{1,2\}$,
we define a cochain $\hat{f}\in C^{n}(\A_{1}\oplus\A_{2})$ by
\begin{eqnarray*}
\hat{f}:=
\left\{
\begin{array}{ll}
f & \text{on $\A_{i(1)}\ot\A_{i(2)}\ot...\ot\A_{i(n)}$,} \\
0  & \text{all other cases.}
\end{array}
\right.
\end{eqnarray*}
We call the cochain $\hat{f}$ a horizontal
lift of $f$, or simply, lift.
For instance, the lifts
of $\a:\A_{1}\ot\A_{1}\to\A_{1}$,
$\b:\A_{1}\ot\A_{2}\to\A_{2}$ and
$\gamma:\A_{2}\ot\A_{1}\to\A_{2}$
are defined by, respectively,
\begin{eqnarray}
\label{exabg1}\hat{\a}((a,x),(b,y))&=&(\a(a,b),0),\\
\label{exabg2}\hat{\b}((a,x),(b,y))&=&(0,\b(a,y)),\\
\label{exabg3}\hat{\gamma}((a,x),(b,y))&=&(0,\gamma(x,b)).
\end{eqnarray}
Let $H:\A_{2}\to\A_{1}$ (resp. $H:\A_{1}\to\A_{2}$)
be a 1 cochain. The lift is defined by
$$
\what{H}(a,x)=(H(x),0) \ \ \ (\text{resp. $\what{H}(a,x)=(0,H(a))$}).
$$
For any $(a,x)\in\A_{1}\oplus\A_{2}$,
we have $\what{H}\what{H}(a,x)=\what{H}(H(x),0)=(0,0)$.
\begin{lemma}\label{keylamma}
$\what{H}\what{H}=0$.
\end{lemma}
This lemma will be used in Section 4.
\medskip\\
\indent
We denote by $\A^{l,k}$
the direct sum of all $l+k$-tensor powers of $\A_{1}$ and $\A_{2}$,
where $l$ (resp. $k$) is the number of $\A_{1}$ (resp. $\A_{2}$).
For instance,
$$
\A^{1,2}:=
(\A_{1}\ot\A_{2}\ot\A_{2})\oplus
(\A_{2}\ot\A_{1}\ot\A_{2})\oplus
(\A_{2}\ot\A_{2}\ot\A_{1}).
$$
The tensor space $(\A_{1}\oplus\A_{2})^{\ot n}$
is expanded into the direct sum of $\A^{l,k}$, $l+k=n$.
For instance,
$$
(\A_{1}\oplus\A_{2})^{\ot 2}=
\A^{2,0}\oplus\A^{1,1}\oplus\A^{0,2}.
$$
We consider the space of cochains,
$C^{n}(\A_{1}\oplus\A_{2}):=
\Hom_{\mathbb{K}}((\A_{1}\oplus\A_{2})^{\ot n},\A_{1}\oplus\A_{2})$.
By the standard properties of $\Hom$-functor,
we have
\begin{equation}\label{dectotal}
C^{n}(\A_{1}\oplus\A_{2})\cong
\sum_{l+k=n}C^{n}(\A^{l,k},\A_{1})\oplus
\sum_{l+k=n}C^{n}(\A^{l,k},\A_{2}),
\end{equation}
where the isomorphism is the horizontal lift.\\
\indent
Let $f$ be a $n$-cochain in $C^{n}(\A_{1}\oplus\A_{2})$.
We say the {\em bidegree} of $f$ is $k|l$,
if $f$ is an element in $C^{n}(\A^{l,k-1},\A_{1})$ or
in $C^{n}(\A^{l-1,k},\A_{2})$, where $n=l+k-1$.
We denote the bidegree of $f$ by $||f||=k|l$.
In general, cochains do not have bidegree.
We call a cochain $f$ a {\em homogeneous cochain},
if $f$ has the bidegree.\\
\indent
We have $k+l\ge 2$, because $n\ge 1$.
Thus there are no cochains of
bidegree $0|0$ or $1|0$ or $0|1$.
If the dimension of $\A_{1}$ is finite and $\A_{2}=\A^{*}_{1}$
is the dual space of $\A_{1}$, then a $k|l$-cochain
is identified with
an element in $\A_{1}^{\ot k}\ot\A_{1}^{*\ot l}$.
Hence the definition above is compatible with the classical one.
For instance, the lift of $H:\A_{2}\to\A_{1}$,
$\what{H}\in C^{1}(\A_{1}\oplus\A_{2})$, has the bidegree $2|0$.
We recall
$\hat{\a},\hat{\b},\hat{\gamma}\in C^{2}(\A_{1}\oplus\A_{2})$
in (\ref{exabg1}), (\ref{exabg2}) and (\ref{exabg3}).
One can easily see 
$||\hat{\a}||=||\hat{\b}||=||\hat{\gamma}||=1|2$.
Thus the sum of $\hat{\a}$, $\hat{\b}$ and $\hat{\gamma}$,
\begin{equation}\label{defmu}
\hat{\mu}:=\hat{\a}+\hat{\b}+\hat{\gamma}
\end{equation}
is a homogeneous cochain with bidegree $1|2$.
The cochain $\hat{\mu}$ is a multiplication of
semidirect product type,
$$
\hat{\mu}((a,x),(b,y))=(\a(a,b),\b(a,y)+\gamma(x,b)),
$$
where $(a,x),(b,y)\in\T$.
Remark that $\hat{\mu}$ is not lift (there is no $\mu$),
however, we will use this symbol,
because $\hat{\mu}$ is an interesting homogeneous cochain.\\
\indent
Clearly, the lemma below holds.
\begin{lemma}
Let $f\in C^{n}(\A_{1}\oplus\A_{2})$ be a cochain.
The bidegree of $f$ is $k|l$ if and only if
the following 4 conditions hold.
\begin{description}
	\item[(deg1)] $k+l-1=n$.
	\item[(deg2-1)] If $\mathbf{x}$ is an element in $\A^{l,k-1}$,
	then $f(\mb{x})$ is in $\A_{1}$.
	\item[(deg2-2)] If $\mathbf{x}$ is an element in $\A^{l-1,k}$,
	then $f(\mathbf{x})$ is in $\A_{2}$.
	\item[(deg3)] All the other cases, $f(\mb{x})=0$.
\end{description}
\end{lemma}
\begin{lemma}
If $||f||=k|0$ (resp. $0|k$) and $||g||=l|0$ (resp. $0|l$),
then $\{f,g\}=0$, or simply,
$$
\{(k|0),(l|0)\}=\{(0|k),(0|l)\}=0.
$$
\end{lemma}
\begin{proof}
Assume that $||f||=k|0$ and $||g||=l|0$.
Then $f$ and $g$ are both horizontal lifts of
cochains in $C^{*}(\A_{2},\A_{1})$.
Thus, from the definition of lift, we have $f\c_{i}g=g\c_{j}f=0$
for any $i,j$.
\end{proof}
\begin{lemma}
Let $f\in C^{|f|}(\A_{1}\oplus\A_{2})$
and $g\in C^{|g|}(\A_{1}\oplus\A_{2})$ homogeneous cochains
with bidegrees
$k_{f}|l_{f}$ and $k_{g}|l_{g}$, respectively,
where $|f|$ and $|g|$ are usual degrees of cochains $f$ and $g$.
The composition $f\c_{i}g$
is again a homogeneous cochain, and the bidegree
is $k_{f}+k_{g}-1|l_{f}+l_{g}-1$.
\end{lemma}
\begin{proof}
We show the conditions (deg1)-(deg3).
The condition (deg1) holds, because
$k_{f}+k_{g}-1+l_{f}+l_{g}-1=|f|+|g|=|f\c_{i}g|+1$.
We show the condition (deg2).
Take an element $\mb{x}\ot\mb{y}\ot\mb{z}$ in
$\A^{l_{f}+l_{g}-1,k_{f}+k_{g}-2}$.
We consider
$$
f\c_{i}g(\mathbf{x},\mathbf{y},\mathbf{z})=
f(\mathbf{x},g(\mathbf{y}),\mathbf{z}).
\eqno{(\star)}
$$
If ($\star$) is zero, then it is in $\A_{1}$.
Namely (deg2-1) is satisfied.
So we assume $(\star)\neq 0$.
We consider the case of $g(\mathbf{y})\in\A_{1}$.
In this case, $\mb{y}$ is in $\A^{l_{g},k_{g}-1}$.
and $\mathbf{x}\ot\mathbf{z}$
is in $\A^{l_{f}-1,k_{f}-1}$.
Thus $\mathbf{x}\ot g(\mathbf{y})\ot\mathbf{z}$
is an element in $\A^{l_{f},k_{f}-1}$ which implies
$f(\mathbf{x}\ot g(\mathbf{y})\ot\mathbf{z})\in\A_{1}$.
When the case of $g(\mathbf{y})\in\A_{2}$,
$\mathbf{y}$ is in $\A^{l_{g}-1,k_{g}}$
and $\mathbf{x}\ot\mathbf{z}$ is in $\A^{l_{f},k_{f}-2}$.
Thus $\mathbf{x}\ot g(\mathbf{y})\ot\mathbf{z}$
is an element in $\A^{l_{f},k_{f}-1}$
which gives $f(\mathbf{x}\ot g(\mathbf{y})\ot\mathbf{z})\in\A_{1}$.
Similar way, when $\mathbf{x}\ot\mathbf{y}\ot\mathbf{z}$
is an element in $\A^{l_{f}+l_{g}-2,k_{f}+k_{g}-1}$,
the condition holds.
We show (deg3). If $\mb{x}\ot\mb{y}\ot\mb{z}$ is an element
in $\A^{l_{f}+l_{g}-1+i,k_{f}+k_{g}-2-i}$ and $g(\mathbf{y})\neq 0$,
then $\mathbf{x}\ot g(\mathbf{y})\ot\mathbf{z}$
is in $\A^{l_{f}+i,k_{f}-1-i}$.
When $i\neq 0$, from the assumption,
$f(\mathbf{x}\ot g(\mathbf{y})\ot\mathbf{z})=0$.
The proof is completed.
\end{proof}
\begin{proposition}
If $||f||=k_{f}|l_{f}$ and $||g||=k_{g}|l_{g}$,
then the Gerstenhaber bracket $\{f,g\}$ has the bidegree
$k_{f}+k_{g}-1|l_{f}+l_{g}-1$.
\end{proposition}
\begin{proof}
Straightforward.
\end{proof}
\noindent
\textbf{Remark}.
Given a bidegree $k+1|l+1$-cochain $f$,
we define $bideg(f):=k|l$.
If $bideg(f)=k|l$ and $bideg(g)=m|n$,
then $bideg(\{f,g\})=bideg(f)+bideg(g)=k+m|l+n$.
Thus the bidegree, $bideg$, of Gerstenhaber bracket is $0|0$.

\section{Main objects.}

\textbf{Notations}.
Let $\A_{1}$ and $\A_{2}$ be vector spaces.
We denote any elements of $\A_{1}$
by $a,b,c,...$ and denote any elements of $\A_{2}$
by $x,y,z,...$.
We sometimes use an identification $(a,x)\cong a+x$,
where $(a,x)\in\A_{1}\oplus\A_{2}$.

\subsection{Twilled algebras.}

\subsubsection{Structures.}

Let $\T$ be an associative algebra equipped with
an associative structure $\theta$.
We assume a decomposition of $\T$ into two subspaces,
$\T=\A_{1}\oplus\A_{2}$.
The associative structure defines
an associative multiplication
by $\theta((a,x),(b,y)):=(a,x)*(b,y)$,
for any $(a,x),(b,y)\in\T$.
\begin{definition}
(\cite{CGM})
The triple $(\T,\A_{1},\A_{2})$, or simply $\T$, is called
an associative \textbf{twilled algebra},
if $\A_{1}$ and $\A_{2}$ are subalgebras of $\T$.
We sometimes denote a twilled algebra $\T$ by $\A_{1}\Join\A_{2}$.
\end{definition}
One can easily check that if $\A_{1}\Join\A_{2}$ is a twilled algebra,
then $\A_{1}$ (resp. $\A_{2}$) is an $\A_{2}$-bimodule
(resp. $\A_{1}$-bimodule). These bimodule structures are defined by
the following decomposition of associative multiplication of $\T$.
For any $a\in\A_{1}$ and $x\in\A_{2}$, the multiplications $a*x$
and $x*a$ are decomposed into 4 multiplications,
$$
a*x=(a*_{2}x, a*_{1}x), \ \
x*a=(x*_{2}a, x*_{1}a),
$$
where $a*_{2}x$ and $x*_{2}a$ are $\A_{1}$-components
of $a*x$ and $x*a$ respectively, and similar way,
$a*_{1}x$ and $x*_{1}a$ are $\A_{2}$-components.
One can easily check that
the multiplication $*_{1}$ (resp. $*_{2}$) is the bimodule action of
$\A_{1}$ to $\A_{2}$ (resp. $\A_{2}$ to $\A_{1}$).\\
\indent
In general, the associative multiplication
of $\A_{1}\Join\A_{2}$ has the form,
$$
(a,x)*(b,y)=(a*b+a*_{2}y+x*_{2}b,a*_{1}y+x*_{1}b+x*y).
$$
The total multiplication, $*$, is decomposed into
two ``associative" multiplications of semidirect product,
\begin{eqnarray*}
(a,x)*_{1}(b,y)&:=&(a*_{1}b,a*_{1}y+x*_{1}b),\\
(a,x)*_{2}(b,y)&:=&(a*_{2}y+x*_{2}b,x*_{2}y),
\end{eqnarray*}
where we put $a*_{1}b:=a*b$ and $x*_{2}y:=x*y$.
Hence the structure $\theta$ is also
decomposed into two associative structures,
$$
\theta=\hat{\mu}_{1}+\hat{\mu}_{2},
$$
where $\hat{\mu}_{i}$ is the structure associated with
the multiplication $*_{i}$ for $i=1,2$.
Recall (\ref{defmu}).
The cochains $\hat{\mu}_{1}$ and $\hat{\mu}_{2}$
have the bidegrees $1|2$ and $2|1$ respectively.
Under the assumption,
the decomposition of $\theta$ is unique, i.e.,
if $\theta$ is decomposed into two substructures of
bidegrees $1|2$ and $2|1$, then such substructures
are uniquely determined.
\begin{lemma}
The associativity of $\theta$ ($\{\theta,\theta\}=0$)
is equivalent with the compatibility conditions,
\begin{eqnarray}
\label{good1}\frac{1}{2}\{\hat{\mu}_{1},\hat{\mu}_{1}\}&=&0,\\
\label{good2}\{\hat{\mu}_{1},\hat{\mu}_{2}\}&=&0,\\
\label{good3}\frac{1}{2}\{\hat{\mu}_{2},\hat{\mu}_{2}\}&=&0.
\end{eqnarray}
\end{lemma}
\begin{proof}
We will show a more generalized result in Lemma \ref{5decomp} below.
\end{proof}
\subsubsection{The case of subalgebras in duality}
Given an arbitrary associative algebra $\A$,
we have a Lie algebra by the commutator, $[a,b]:=ab-ba$ on $\A$.
The induced Lie algebra is denoted by $L(\A)$.
The correspondence
$L:\A\to L(\A)$ is a functor (sometimes called a Liezation)
from the usual category of associative algebras
to the one of Lie algebras.
\medskip\\
\indent
In this short section, we assume that
$\A_{1}=:\A$ is a finite dimensional vector space and $\A_{2}$
is the dual space. In this case, $\T=\A\oplus\A^{*}$
has a nondegenerate symmetric bilinear form, $(-|-)$,
where $(\A|\A^{*})=(\A^{*}|\A)$ is the dual pairing
and $(\A|\A)=(\A^{*}|\A^{*})=0$.
We set a natural assumption, namely,
the bilinear form is invariant (or associative)
with respect to the associative multiplication of $\T$,
explicitly,
$$
(t_{1}*t_{2}|t_{3})=(t_{1}|t_{2}*t_{3})
$$
for any $t_{1},t_{2},t_{3}\in\T$. Such a twilled
algebra is called an invariant twilled algebra.\\
\indent
If $\T$ is an invariant twilled algebra,
then the triple $(L(\T),L(\A),L(\A^{*}))$ is a Manin triple.
It is a twilled Lie algebra with an invariant
pseudo-Euclidean metric decomposed into
two maximally isotropic subalgebras.
In general, a pair of Lie algebras $(\g_{1},\g_{2})$
becomes a Lie bialgebra if and only if
a triple of Lie algebras
$(\g_{1}\Join\g_{2},\g_{1},\g_{2})$ is
a Manin triple. In this times,
the total space $\g_{1}\Join\g_{2}$
is called a Drinfeld double.
Thus the pair $(L(\A),L(\A^{*}))$ becomes a Lie bialgebra
and $L(\A)\Join L(\A^{*})$ is a Drinfeld double.
If $\T$ is a quasi-twilled algebra
in Definition \ref{defqtwi} below,
then the cocycle term $\phi_{1}$ (or $\phi_{2}$)
is a cyclic cocycle, i.e., for any $a,b,c\in\A$,
$$
\phi_{1}(a,b)(c)=\phi_{1}(b,c)(a)=\phi_{1}(c,a)(b).
$$
This fact is directly checked by the invariancy.
And the commutator, $\Phi_{1}(a,b):=\phi_{1}(a,b)-\phi_{1}(b,a)$,
is identified with a skew symmetric 3-tensor in $\bigwedge^{3}\A^{*}$.
This implies that if $\A\oplus\A^{*}$
is a quasi-twilled algebra, then $L(\T)$
is the double of quasi-Lie bialgebra $(L(\A),L(\A)^{*})$
(see \cite{D},\cite{Kos1} for quasi-Lie bialgebras).
\medskip\\
\indent
The dual map of an associative multiplication
on $\T$ becomes a coassociative multiplication $\T\to\T\ot\T$.
Here $\T$ and $\T\ot\T$ are identified with
$\T^{*}$ and $(\T\ot\T)^{*}$ by the bilinear form.
Since $\hat{\mu}_{i}$ is associative,
the dual map of $\hat{\mu}_{i}$
becomes a coassociative multiplication,
$\Delta_{\hat{\mu}_{i}}:\T\to\T\ot\T$, $i=1,2$.
We rewrite the conditions (\ref{good1}), (\ref{good2}) and (\ref{good3})
by the comultiplications.
(\ref{good1}) and (\ref{good3}) are equivalent with
coassociativity of $\Delta_{\hat{\mu}_{i}}$, $i=1,2$, respectively.
So we consider (\ref{good2}).
We define a $(\T,\hat{\mu}_{1})$-bimodule structure on $\T\ot\T$ by
$t\cdot(\T\ot\T):=(t*_{1}\T)\ot\T$ and 
$(\T\ot\T)\cdot{t}:=\T\ot(\T*_{1}t)$
where $t\in\T$ and
$*_{1}$ is the associative multiplication of $\hat{\mu}_{1}$.
For any $s,t,u,v\in\T$, we have
$$
(\Delta_{\hat{\mu}_{2}}(s*_{1}t)|u\ot v)=(s*_{1}t|u*_{2}v),
$$
where the pairing $(-|-)$ is extended on $\T\ot\T$ by the rule,
$$
(s\ot t|u\ot v):=(s|v)(t|u).
$$
The invariancy holds with respect to $\hat{\mu}_{i}$, $i=1,2$,
for instance,
$$
(a*_{1}x|b)=(a*x|b)=(a|x*b)=(a|x*_{1}b),
$$
where $(\A|\A)=0$ is used.
From the invariancy,
we have $(s*_{1}t|u*_{2}v)=(s|t*_{1}(u*_{2}v))$.
By (\ref{good2}), we have
$t*_{1}(u*_{2}v)=(t*_{2}u)*_{1}v+(t*_{1}u)*_{2}v-t*_{2}(u*_{1}v)$.
Thus (\ref{good2}) is equivalent with the condition,
\begin{multline}\label{aubku1}
(\Delta_{\hat{\mu}_{2}}(s*_{1}t)|u\ot v)=(s|t*_{1}(u*_{2}v))=\\
=(s|(t*_{2}u)*_{1}v)+(s|(t*_{1}u)*_{2}v)-(s|t*_{2}(u*_{1}v)).
\end{multline}
The first term of the right-hand side of (\ref{aubku1}) is
$$
(s|(t*_{2}u)*_{1}v)=(v*_{1}s|t*_{2}u)=
(u*_{2}(v*_{1}s)|t)=
(u\ot(v*_{1}s)|\Delta_{\hat{\mu}_{2}}(t)).
$$
We put $\Delta_{\hat{\mu}_{2}}(t)=\sum t_{1}\ot{t_{2}}$.
Then we have
$$
(u\ot(v*_{1}s)|\Delta_{\hat{\mu}_{2}}(t))
=\sum(u|t_{2})(v*_{1}s|t_{1})=\sum(u|t_{2})(v|s*_{1}t_{1})
=(u\ot v|s\cdot \Delta_{\hat{\mu}_{2}}(t)).
\eqno{(A)}
$$
And the second and third terms of the right-hand side of
(\ref{aubku1}) are
$$
(s|(t*_{1}u)*_{2}v)-(s|t*_{2}(u*_{1}v))=
(\Delta_{\hat{\mu}_{2}}(s)|(t*_{1}u)\ot v)-(s*_{2}t|u*_{1}v).
$$
We put $\Delta_{\hat{\mu}_{2}}(s)=\sum s_{1}\ot s_{2}$.
Then we have
$$
(\Delta_{\hat{\mu}_{2}}(s)|(t*_{1}u)\ot v)
=\sum(s_{1}|v)(s_{2}|t*_{1}u)=
\sum(s_{1}|v)(s_{2}*_{1}t|u)=
(\Delta_{\hat{\mu}_{2}}(s)\cdot{t}|u\ot v).
\eqno{(B)}
$$
and
$$
(s*_{2}t|u*_{1}v)=
(\Delta_{\hat{\mu}_{1}}(s*_{2}t)|u\ot v)=
(\Delta_{\hat{\mu}_{1}}\c\hat{\mu}_{2}(s,t)|u\ot v).
\eqno{(C)}
$$
From (A),(B) and (C), we obtain a compatibility condition,
\begin{equation}\label{mulcocomp}
(\Delta_{\hat{\mu}_{2}}(s*_{1}t)|u\ot v)=
(s\cdot \Delta_{\hat{\mu}_{2}}(t)|u\ot v)+
(\Delta_{\hat{\mu}_{2}}(s)\cdot{t}|u\ot v)
-(\Delta_{\hat{\mu}_{1}}\c\hat{\mu}_{2}(s,t)|u\ot v).
\end{equation}
Since $\T\ot\T$ is a $(\T,\hat{\mu}_{1})$-bimodule,
we have a Hochschild complex
$(C^{*}(\T,\T\ot\T),D_{\hat{\mu}_{1}})$,
where $D_{\hat{\mu}_{1}}$ is a Hochschild coboundary map.
The condition (\ref{mulcocomp}) is equivalent with (\ref{pdcom2}) below.
Under the assumptions of this section,
the identity (\ref{good2}) $\{\hat{\mu}_{1},\hat{\mu}_{2}\}=0$ is
equivalent with
\begin{equation}\label{pdcom2}
D_{\hat{\mu}_{1}}\Delta_{\hat{\mu}_{2}}
-\Delta_{\hat{\mu}_{1}}\c\hat{\mu}_{2}=0.
\end{equation}
Since $\{\hat{\mu}_{2},\hat{\mu}_{1}\}=0$, we have
$D_{\hat{\mu}_{2}}\Delta_{\hat{\mu}_{1}}
-\Delta_{\hat{\mu}_{2}}\c\hat{\mu}_{1}=0$.
One can easily show that $D_{\hat{\mu}_{i}}\Delta_{\hat{\mu}_{i}}
-\Delta_{\hat{\mu}_{i}}\c\hat{\mu}_{i}=0$
holds for $i=1,2$.
Thus we have $D_{\theta}\Delta_{\theta}-\Delta_{\theta}\c\theta=0$.
From (\ref{pdcom2}) we have
$D_{\hat{\mu}_{1}}(\Delta_{1}\c\hat{\mu}_{2})=0$.
By direct computation, one can show that
if $\A$ is unital (i.e. $1*_{1}\A=\A*_{1}1$), then
$D_{\hat{\mu}_{1}}(\Delta_{1}\c\hat{\mu}_{2})=0$
implies (\ref{pdcom2}).
\medskip\\
\indent
It is obvious that $\A$ is a sub-coalgebra of
$(\T,\Delta_{\hat{\mu}_{2}})$.
Since $\hat{\mu}_{2}$ is zero on $\A\ot\A$,
$\Delta_{\hat{\mu}_{2}}$ is a derivation on $\A$, i.e.,
for any $a,b\in\A$,
$$
\Delta_{\hat{\mu}_{2}}(a*_{1}b)=
\Delta_{\hat{\mu}_{2}}(a)\cdot b+a\cdot\Delta_{\hat{\mu}_{2}}(b).
$$
An associative and coassociative algebra
$(\I,*,\delta)$ is called an infinitesimal bialgebra (\cite{JR}),
if $\delta(a*b)=a\cdot\delta(b)+\delta(b)\cdot a$ for any $a,b\in\I$.
Thus the triple $(\A,*_{1},\Delta_{\hat{\mu}_{2}})$
is an infinitesimal bialgebra.
We consider the converse.
Given an infinitesimal bialgebra $(\I,*,\delta)$,
the multiplications $*$ and $\delta$ are extended
on $\I\oplus\I^{*}$ by adjoint actions.
However the compatibility condition (\ref{pdcom2})
is not satisfied in general. This implies that the Liezation of
an infinitesimal bialgebra is not a Lie bialgebra in general.
For this problem, see the detailed study Aguiar \cite{Ag2}.

\subsubsection{Induced dg-Lie algebras.}

This short section is the heart of this article.
The meaning of twilled algebra is given by
the proposition below.
From the associative condition (\ref{good1}),
$(C^{*}(\T),d_{\hat{\mu}_{1}}(-):=\{\hat{\mu}_{1},-\})$
becomes a dg-Lie algebra.
The graded space $C^{*}(\A_{2},\A_{1})$ is identified
with an abelian subalgebra of the dg-Lie algebra,
via the horizontal lift. 
One can easily check that the derived bracket
$$
[f,g]_{\hat{\mu}_{1}}:=(-1)^{|f|-1}\{\{\hat{\mu}_{1},f\},g\}
$$
is closed on $C^{*}(\A_{2},\A_{1})$.
From Lemma \ref{balma}, $C^{*}(\A_{2},\A_{1})$
becomes a graded Lie algebra.
Further, by (\ref{good2}) and (\ref{good3}),
$d_{\hat{\mu}_{2}}:=\{\hat{\mu}_{2},\}$
becomes a square zero derivation on
the induced graded Lie algebra $C^{*}(\A_{2},\A_{1})$.
\begin{proposition}\label{dglie}
If $\T=\A_{1}\Join\A_{2}$ is a twilled algebra, then
$C^{*}(\A_{2},\A_{1})$ has a dg-Lie algebra structure.
The degree of dg-Lie algebra structure is the same as
the usual degree of cochains.
\end{proposition}
\begin{proof}
We show only a derivation property of $d_{\hat{\mu}_{2}}$.
Since $\hat{\mu}_{2}$ is an associative structure,
$d_{\hat{\mu}_{2}}$ is square zero.
For any cochains $f,g\in C^{*}(\A_{2},\A_{1})$,
we have
\begin{eqnarray*}
d_{\hat{\mu}_{2}}[f,g]_{\hat{\mu}_{1}}&:=&(-1)^{|f|-1}
\{\hat{\mu}_{2},\{\{\hat{\mu}_{1},f\},g\}\} \\
&=&(-1)^{|f|-1}\{\{\hat{\mu}_{2},\{\hat{\mu}_{1},f\}\},g\}
-\{\{\hat{\mu}_{1},f\},\{\hat{\mu}_{2},g\}\} \\
&=&(-1)^{|f|}\{\{\hat{\mu}_{1},\{\hat{\mu}_{2},f\}\},g\}
-\{\{\hat{\mu}_{1},f\},\{\hat{\mu}_{2},g\}\} \\
&=&[d_{\hat{\mu}_{2}}f,g]_{\hat{\mu}_{1}}+
(-1)^{|f|}[f,d_{\hat{\mu}_{2}}g]_{\hat{\mu}_{1}}.
\end{eqnarray*}
From Lemma \ref{balma},
the derived degree is given by
$deg_{d_{\hat{\mu}_{1}}}(f)=deg(f)+1=|f|$,
where $deg(f)=|f|-1$ is the degree of the canonical
dg-Lie algebra $(C^{*}(\T),d_{\hat{\mu}_{1}})$
(recall Section 2.1).
Thus $d_{\hat{\mu}_{2}}$ satisfies the defining condition
(\ref{defdgl3}) of dg-Lie algebra.
\end{proof}
When we recall deformation theory,
it is natural to ask: What is a solution of
Maurer-Cartan equation
in the dg-Lie algebra ? We will solve this question
in Section 5.
\subsubsection{Examples}
\begin{example}\label{triext}
(trivial extensions, semidirect product algebras.)
Let $\A$ be an associative algebra and let $M$ an $\A$-bimodule.
The trivial extension $\A\ltimes M$
is a twilled algebra of $\A=\A_{1}$ and $M=\A_{2}$,
where the structure $\hat{\mu}_{2}$ is trivial and $\hat{\mu}_{1}$
is defined by, for any $(a,m),(b,n)\in\A\oplus M$,
$$
\hat{\mu}_{1}((a,m),(b,n)):=
(a,m)*(b,n):=(ab,a\cdot n+m\cdot b),
$$
where $\cdot$ is the bimodule action of $\A$ on $M$.
\end{example}
A direct product algebra $\A\times\A$ is a twilled algebra.
The following example is considered as
a $q$-analogue of trivial extensions.
\begin{example}\label{nonabext}
(q-trivial extensions.)
Let $\A$ be an associative algebra.
Define a multiplication on $\A\oplus\A$ by
$$
(a,x)*_{q}(b,y):=(ab,ay+xb+qxy),
$$
where $q\in\mathbb{K}$.
Then $(\A\oplus\A,*_{q})$ becomes a twilled algebra.
We denote the twilled algebra by $\A\Join_{q}\A$.
\end{example}
If $(\T,\theta)$ is an associative
algebra, then $C^{*}(\T)$ becomes an associative
algebra by a cup product,
$f\vee_{\theta}g:=\theta(f,g)$, $f,g\in C^{*}(\T)$.
\begin{example}
If $\T=\A_{1}\Join\A_{2}$ is a twilled algebra, then
$$
C^{*}(\T)=C^{*}(\T,\A_{1}\Join\A_{2})\cong
C^{*}(\T,\A_{1})\Join C^{*}(\T,\A_{2})
$$
is a twilled algebra, because
the cup product is decomposed into
$\vee_{\theta}=\vee_{\hat{\mu}_{1}}+\vee_{\hat{\mu}_{2}}$.
\end{example}

\subsection{Proto-, Quasi-twilled algebras.}

A {\em quasi-Lie} bialgebra is known as
a classical limit of a quasi-Hopf algebra.
The notion of quasi-Lie bialgebra is generalized
to {\em proto-Lie} bialgebras (see \cite{Kos1}).
The latter is more complicated object
than quasi-Lie bialgebras.
The proto-Lie bialgebras
provide a general framework of
quantum-classical correspondence.
In this section, we will study associative analogues
of proto-, quasi-Lie bialgebras.
\begin{definition}
Let $(\T,\theta)$ be an associative algebra
decomposed into two subspaces,
$\T=\A_{1}\oplus\A_{2}$.
Here $\A_{1}$ and $\A_{2}$ are not
necessarily subalgebras.
We call the triple $(\T,\A_{1},\A_{2})$
a \textbf{proto-twilled algebra}.
\end{definition}

\begin{lemma}\label{dec4str}
Let $\theta$ be an arbitrary 2-cochain in $C^{2}(\T)$.
Then $\theta$ is uniquely decomposed into 4 homogeneous cochains
of bidegrees $0|3$, $1|2$, $2|1$ and $3|0$,
\begin{equation*}
\theta=\hat{\phi}_{1}+\hat{\mu}_{1}+\hat{\mu}_{2}+\hat{\phi}_{2}.
\end{equation*}
\end{lemma}
\begin{proof}
Recall the decomposition (\ref{dectotal}).
The space of $2$-cochains
$C^{2}(\T)$ is decomposed into 4 subspaces,
$$
C^{2}(\T)=(0|3)\oplus(1|2)\oplus(2|1)\oplus(3|0),
$$
where $(i|j)$ is the space of bidegree $i|j$-cochains, $i,j=0,1,2,3$.
The decomposition is essentially unique.
Thus $\theta$ is uniquely decomposed into
homogeneous cochains of bidegrees $0|3$, $1|2$, $2|1$ and $3|0$.
The 4 substructures $\hat{\phi}_{1}$,
$\hat{\mu}_{1}$, $\hat{\mu}_{2}$ and $\hat{\phi}_{2}$
in the lemma are given as the homogeneous cochains.
The proof is completed.
\end{proof}
The multiplication $(a,x)*(b,y):=\theta((a,x),(b,y))$ of $\T$
is uniquely decomposed by the canonical projections
$\T\to\A_{1}$ and $\T\to\A_{2}$ into the 8 multiplications,
\begin{eqnarray*}
a*b&=&(a*_{1}b,a*_{2}b),\\
a*y&=&(a*_{2}y,a*_{1}y),\\
x*b&=&(x*_{2}b,x*_{1}b),\\
x*y&=&(x*_{1}y,x*_{2}y).
\end{eqnarray*}
We put bidegrees on the 4 cochains,
$\|\hat{\phi}_{1}\|:=0|3$,
$\|\hat{\mu}_{1}\|:=1|2$,
$\|\hat{\mu}_{2}\|:=2|1$ and
$\|\hat{\phi}_{2}\|:=3|0$.
Then we obtain
\begin{eqnarray*}
\hat{\phi}_{1}((a,x),(b,y))&=&(0,a*_{2}b),\\
\hat{\mu}_{1}((a,x),(b,y))&=&(a*_{1}b,a*_{1}y+x*_{1}b),\\
\hat{\mu}_{2}((a,x),(b,y))&=&(a*_{2}y+x*_{2}b,x*_{2}y),\\
\hat{\phi}_{2}((a,x),(b,y))&=&(x*_{1}y,0).
\end{eqnarray*}
Remark that $\hat{\phi}_{1}$ and $\hat{\phi}_{2}$ are lifted cochains
of $\phi_{1}(a,b):=a*_{2}b$ and $\phi_{2}(x,y):=x*_{1}y$.
\begin{lemma}\label{5decomp}
The Maurer-Cartan condition $\{\theta,\theta\}=0$
is equivalent with the following 5 conditions.
\begin{eqnarray}
\label{5cond4}\{\hat{\mu}_{1},\hat{\phi}_{1}\}&=&0,\\
\label{5cond1}\frac{1}{2}\{\hat{\mu}_{1},\hat{\mu}_{1}\}+
\{\hat{\mu}_{2},\hat{\phi}_{1}\}&=&0,\\
\label{5cond2}\{\hat{\mu}_{1},\hat{\mu}_{2}\}+
\{\hat{\phi}_{1},\hat{\phi}_{2}\}&=&0,\\
\label{5cond3}\frac{1}{2}\{\hat{\mu}_{2},\hat{\mu}_{2}\}+
\{\hat{\mu}_{1},\hat{\phi}_{2}\}&=&0,\\
\label{5cond5}\{\hat{\mu}_{2},\hat{\phi}_{2}\}&=&0.
\end{eqnarray}
\end{lemma}
\begin{proof}
From the 5 conditions, one can directly check
the Maurer-Cartan condition of $\theta$.
We show the converse.
The bidegrees of $\hat{\phi}_{1}$,
$\hat{\mu}_{1}$, $\hat{\mu}_{2}$
and $\hat{\phi}_{2}$
are $0|3$, $1|2$, $2|1$ and $3|0$, respectively.
If $\{\theta,\theta\}=0$, then
\begin{multline*}
\{\hat{\mu}_{1},\hat{\mu}_{1}\}+2\{\hat{\mu}_{2},\hat{\phi}_{1}\}+
2\{\hat{\mu}_{1},\hat{\mu}_{2}\}+2\{\hat{\phi}_{1},\hat{\phi}_{2}\}+
\{\hat{\mu}_{2},\hat{\mu}_{2}\}+2\{\hat{\mu}_{1},\hat{\phi}_{2}\}+\\
2\{\hat{\mu}_{1},\hat{\phi}_{1}\}+
2\{\hat{\mu}_{2},\hat{\phi}_{2}\}=0.
\end{multline*}
The first two terms have $1|3$-bidegree,
the second two terms have $2|2$-bidegree,
the third two terms have $3|1$-bidegree
and the last two terms have $0|4$ and $4|0$ respectively.
Thus we have
$\{\hat{\mu}_{1},\hat{\mu}_{1}\}+2\{\hat{\mu}_{2},\hat{\phi}_{1}\}=0$
for $1|3$-bidegree, and this is (\ref{5cond1}).
Similarly, we obtain (\ref{5cond4})-(\ref{5cond5}).
\end{proof}
\begin{definition}\label{defqtwi}
Let $\T=\A_{1}\oplus\A_{2}$
be a proto-twilled algebra equipped with the structures
$(\hat{\mu}_{1},\hat{\mu}_{2},\hat{\phi}_{1},\hat{\phi}_{2})$.
We call the triple $(\T,\A_{1},\A_{2})$
a \textbf{quasi-twilled algebra}, if $\phi_{2}=0$,
or equivalently, $\A_{2}$ is a subalgebra.
Since $\A_{1}\oplus\A_{2}=\A_{2}\oplus\A_{1}$,
the definition is adapted in the case of
$\phi_{2}\neq 0$ and $\phi_{1}=0$.
\end{definition}
It is obvious that twilled algebras are special quasi-twilled algebras
of $\phi_{1}=\phi_{2}=0$.
From Lemma \ref{5decomp}, $\theta$ is the structure of
a quasi-twilled algebra of $\phi_{2}=0$ if and only if
\begin{eqnarray}
\label{qcon4}\{\hat{\mu}_{1},\hat{\phi}_{1}\}&=&0, \\
\label{qcon1}\frac{1}{2}\{\hat{\mu}_{1},\hat{\mu}_{1}\}+
\{\hat{\mu}_{2},\hat{\phi}_{1}\}&=&0,\\
\label{qcon2}\{\hat{\mu}_{1},\hat{\mu}_{2}\}&=&0,\\
\label{qcon3}\frac{1}{2}\{\hat{\mu}_{2},\hat{\mu}_{2}\}&=&0.
\end{eqnarray}
In Proposition \ref{dglie},
we saw $C^{*}(\A_{2},\A_{1})$ has a dg-Lie algebra structure.
In the quasi-twilled algebra cases,
from (\ref{qcon3}), $d_{\hat{\mu}_{2}}$ is still a square zero
derivation, but the derived bracket by $\hat{\mu}_{1}$
does not satisfy the graded Jacobi identity in general.
However the Jacobiator
still satisfies a weak Jacobi identity
in the sense of homotopy Lie algebras (\cite{DMZ},\cite{LM}).
The 3-cochain $\frac{1}{2}\{\hat{\mu}_{1},\hat{\mu}_{1}\}$
rises up to the graded Jacobiator
via the derived bracket,
\begin{multline*}
(-1)^{|g|-1}
\frac{1}{2}\{\{\{\{\hat{\mu}_{1},\hat{\mu}_{1}\},f\},g\},h\}=\\
[f,[g,h]_{\hat{\mu}_{1}}]_{\hat{\mu}_{1}}
-[[f,g]_{\hat{\mu}_{1}},h]_{\hat{\mu}_{1}}
-(-1)^{|f||g|}[g,[f,h]_{\hat{\mu}_{1}}]_{\hat{\mu}_{1}}.
\end{multline*}
From (\ref{qcon1}), the Jacobiator is also given by
$-\{\hat{\mu}_{2},\hat{\phi}_{1}\}$.
We define a tri-linear bracket product (homotopy)
on $C^{*}(\A_{2},\A_{1})$ by
$$
[f,g,h]_{\hat{\phi}_{1}}:=(-1)^{|g|-1}
\{\{\{\hat{\phi}_{1},f\},g\},h\}.
$$
Since $C^{*}(\A_{2},\A_{1})$ is abelian with respect to $\{-,-\}$,
the tribracket is skew-symmetric.
We can show that the system,
$(d_{\hat{\mu}_{2}},
[\cdot,\cdot]_{\hat{\mu}_{1}},
[\cdot,\cdot,\cdot]_{\hat{\phi}_{1}})$,
defines a strong homotopy
Lie algebra structure of $l_{n\ge 4}:=0$
on $C^{*}(\A_{2},\A_{1})$.
This assertion will be shown as a corollary
of a more general result in \cite{U2}.
\medskip\\
\indent
The complex plane, $\T:=\mathbb{C}$,
is a quasi-twilled algebra decomposed into
the real part and the imaginary part.
Given a $\mathbb{R}$-algebra $\A$,
the complexification
$\mathbb{C}\ot_{\mathbb{R}}\A=\A\oplus\sqrt{-1}\A$
is a quasi-twilled algebra.
\begin{example}\label{nonabext2}
(Quasi-trivial extension.)
Let $\A$ be an associative algebra.
Define a multiplication on $\A\oplus\A$ by
$$
(a,x)*_{Q}(b,y):=(ab+Qxy,ay+xb),
$$
where $Q\in\mathbb{K}$.
Then $\A\oplus\A$ becomes a quasi-twilled algebra,
where $\phi_{2}(x,y):=Qxy$.
We denote the algebra by $\A\oplus_{Q}\A$.
\end{example}

\section{Twisting by a 1-cochain}

Let $h$ be a 1-cochain in $C^{1}(\T)$.
By analogy with Hamiltonian vector field,
we define an operator by
$X_{h}:=\{\cdot,h\}$,
and by analogy with Hamiltonian flow, we put
$$
exp(X_{h})(\cdot):=1+X_{h}+\frac{1}{2!}X^{2}_{h}
+\frac{1}{3!}X^{3}_{h}+...,
$$
where $X^{2}_{h}:=\{\{\cdot,h\},h\}$
and $X^{n}_{h}$ is defined by the same manner.
Remark that $exp(X_{h})$ is not well-defined in general.
\medskip\\
\indent
Let $(\T=\A_{1}\oplus\A_{2},\theta)$
be a proto-twilled algebra,
and let $\hat{H}\in C^{1}(\T)$ be the lift of a linear map
$H:\A_{2}\to\A_{1}$ (or $H:\A_{1}\to\A_{2}$).
Then $exp(X_{\what{H}})$ is always well-defined
as an operator, because $\what{H}\what{H}=0$
(recall Lemma \ref{keylamma}).
\begin{definition}\label{corodef}
A transformation (\ref{thoverH}) is called
a ``twisting" of $\theta$ by $H$.
\begin{equation}\label{thoverH}
\theta^{H}:=exp(X_{\what{H}})(\theta).
\end{equation}
\end{definition}
It is clear that the result of twisting by $H$
is again a 2-cochain.
We can consider the twisting operations
are special examples of gauge transformations
in deformation theory (see \cite{DMZ}).
The following Lemma \ref{stag1} and Proposition \ref{stag2}
are followed from standard arguments in deformation theory.
\begin{lemma}\label{stag1}
$\theta^{H}=e^{-\what{H}}\theta(e^{\what{H}}\ot e^{\what{H}})$,
where $e^{\pm\what{H}}=1\pm\what{H}$.
\end{lemma}
\begin{proof}
We have $e^{-\what{H}}\theta(e^{\what{H}}\ot e^{\what{H}})=
\theta(e^{\what{H}}\ot e^{\what{H}})-
\what{H}\theta(e^{\what{H}}\ot e^{\what{H}})=$
\begin{multline*}
=\theta+\theta(1\ot\what{H})+\theta(\what{H}\ot 1)
+\theta(\what{H}\ot\what{H})
-\what{H}\theta-\what{H}\theta(1\ot\what{H})
-\what{H}\theta(\what{H}\ot 1)
-\what{H}\theta(\what{H}\ot\what{H})=\\
\theta+\theta(1\ot\what{H})+\theta(\what{H}\ot 1)-\what{H}\theta+
\theta(\what{H}\ot\what{H})-\what{H}\theta(1\ot\what{H})
-\what{H}\theta(\what{H}\ot 1)
-\what{H}\theta(\what{H}\ot\what{H}).
\end{multline*}
Since $\what{H}\what{H}=0$, for any $I\ge 4$,
we have $X_{\what{H}}^{I}(\theta)=0$.
Thus we have
$$
exp(X_{\what{H}})(\theta)=
\theta+\{\theta,\what{H}\}+
\frac{1}{2}\{\{\theta,\what{H}\},\what{H}\}+
\frac{1}{6}\{\{\{\theta,\what{H}\},\what{H}\},\what{H}\}.
$$
One can directly check the three identities below.
\begin{eqnarray*}
\{\theta,\what{H}\}&=&
\theta(\what{H}\ot 1)+\theta(1\ot\what{H})-\what{H}\theta,\\
\frac{1}{2}\{\{\theta,\what{H}\},\what{H}\}&=&
\theta(\what{H}\ot\what{H})
-\what{H}\theta(\what{H}\ot 1)-\what{H}\theta(1\ot\what{H}),\\
\frac{1}{6}\{\{\{\theta,\what{H}\},\what{H}\},\what{H}\}&=&
-\what{H}\theta(\what{H}\ot\what{H}).
\end{eqnarray*}
The proof of the lemma is completed.
\end{proof}
From above lemma, we have
$\{\theta^{H},\theta^{H}\}=e^{-H}\{\theta,\theta\}
(e^{\what{H}}\ot e^{\what{H}}\ot e^{\what{H}})$.
This implies
\begin{proposition}\label{stag2}
The result of twisting
$\theta^{H}$ is an associative structure, i.e.,
$\{\theta^{H},\theta^{H}\}=0$.
\end{proposition}
The following corollary is useful.
\begin{corollary}\label{benri}
The twisting by $H$ induces an algebra isomorphism,
$$
e^{H}:(\T,\theta^{H})\to(\T,\theta).
$$
\end{corollary}

Obviously, $(\T,\theta^{H})$ is also
a proto-twilled algebra.
Thus $\theta^{H}$ is also
decomposed into the unique 4 substructures.
The twisting operations are completely determined by
\begin{theorem}\label{maintheorem}
Assume a decomposition of $\theta$,
$\theta:=\hat{\mu}_{1}+\hat{\mu}_{2}+\hat{\phi}_{1}+\hat{\phi}_{2}$.
The unique 4 substructures of $\theta^{H}$ have the following form:
\begin{eqnarray}
\label{4sth3}\hat{\phi}^{H}_{1}&=&\hat{\phi}_{1},\\
\label{4sth1}\hat{\mu}^{H}_{1}&=&
\hat{\mu}_{1}+\{\hat{\phi}_{1},\what{H}\},\\
\label{4sth2}\hat{\mu}^{H}_{2}&=&\hat{\mu}_{2}+
d_{\hat{\mu}_{1}}\what{H}+
\frac{1}{2}\{\{\hat{\phi}_{1},\what{H}\},\what{H}\},\\
\label{4sth4}
\hat{\phi}^{H}_{2}&=&\hat{\phi}_{2}+d_{\hat{\mu}_{2}}\what{H}+
\frac{1}{2}[\what{H},\what{H}]_{\hat{\mu}_{1}}+
\frac{1}{6}\{\{\{\hat{\phi}_{1},\what{H}\},\what{H}\},\what{H}\},
\end{eqnarray}
where $d_{\hat{\mu}_{i}}(-):=
\{\hat{\mu}_{i},-\}$, $(i=1,2)$ and
$[\what{H},\what{H}]_{\hat{\mu}_{1}}:=
\{\{\hat{\mu}_{1},\what{H}\},\what{H}\}$.
\end{theorem}
\begin{proof}
The first term of $exp(X_{\what{H}})(\theta)$ is $\theta$.
From the bidegree calculus,
we have $\{\hat{\phi}_{2},\what{H}\}=0$,
because $||\hat{\phi}_{2}||=3|0$ and
$||\what{H}||=2|0$.
Thus the second term of $exp(X_{\what{H}})(\theta)$ has the form,
$$
\{\hat{\mu}_{1},\what{H}\}+\{\hat{\mu}_{2},\what{H}\}+
\{\hat{\phi}_{1},\what{H}\}.
$$
We have
$\|\{\hat{\mu}_{1},\what{H}\}\|=2|1$,
$\|\{\hat{\mu}_{2},\what{H}\}\|=3|0$ and
$\|\{\hat{\phi}_{1},\what{H}\}\|=1|2$,
which implies $\{\{\hat{\mu}_{2},\what{H}\},\what{H}\}=0$.
Thus the third term has the form,
$$
\frac{1}{2}(\{\{\hat{\mu}_{1},\what{H}\},\what{H}\}+
\{\{\hat{\phi}_{1},\what{H}\},\what{H}\}).
$$
The bidegrees are
$\|\{\{\hat{\mu}_{1},\what{H}\},\what{H}\}\|=3|0$
and $\|\{\{\hat{\phi}_{1},\what{H}\},\what{H}\}\|=2|1$.
The final term is
$\{\{\{\theta,\what{H}\},\what{H}\},\what{H}\}=
\{\{\{\hat{\phi}_{1},\what{H}\},\what{H}\},\what{H}\}$
which has
the bidegree $3|0$.
Thus the sum of all $3|0$-terms is
$$
\hat{\phi}_{2}+\{\hat{\mu}_{2},\what{H}\}+
\frac{1}{2!}\{\{\hat{\mu}_{1},\what{H}\},\what{H}\}+
\frac{1}{3!}\{\{\{\hat{\phi}_{1},\what{H}\},\what{H}\},
\what{H}\}
$$
which gives (\ref{4sth4}).
In this way, the remaining 3 conditions hold.
\end{proof}

\section{Maurer-Cartan equations}

Let $\T=\A_{1}\oplus\A_{2}$ be a proto-twilled algebra
equipped with an associative structure $\theta$ and let 
$(\hat{\phi}_{1},\hat{\mu}_{1},\hat{\mu}_{2},\hat{\phi}_{2})$
be the unique 4 substructures of $\theta$.
In this section, we discuss various examples of 
twisting operations.

\subsection{The cases of $\phi_{1}=0$ and $\phi_{2}=0$.}

In this case, $\T=\A_{1}\Join\A_{2}$ is a twilled algebra.
However the result of twisting by $H:\A_{2}\to\A_{1}$,
$(\T_{H},\A_{1},\A_{2})$, is a quasi-twilled algebra in general.
The twisted structures have the forms,
\begin{eqnarray*}
\hat{\mu}^{H}_{1}&=&\hat{\mu}_{1},\\
\hat{\mu}^{H}_{2}&=&\hat{\mu}_{2}+d_{\hat{\mu}_{1}}\hat{H},\\
\hat{\phi}^{H}_{2}&=&d_{\hat{\mu}_{2}}\what{H}+
\frac{1}{2}[\what{H},\what{H}]_{\hat{\mu}_{1}}.
\end{eqnarray*}
This $\hat{\phi}^{H}_{2}$ is called a curvature.
The derivation operator $d_{\hat{\mu}_{2}}$ on
the graded Lie algebra $C^{*}(\A_{2},\A_{1})$ is modified by $H$,
$d_{\hat{\mu}^{H}_{2}}(-)=
d_{\hat{\mu}_{2}}(-)+[\hat{H},-]_{\hat{\mu}_{1}}$,
where $d_{\hat{\mu}^{H}_{2}}d_{\hat{\mu}^{H}_{2}}\neq 0$ in general.
By Lemma \ref{5decomp} (\ref{5cond5}),
the cocycle condition of $\phi^{H}_{2}$ still holds,
$$
d_{\hat{\mu}^{H}_{2}}\hat{\phi}^{H}_{2}=0.
$$
This is a kind of Bianchi identity.

\subsubsection{Maurer-Cartan operators.}

In Proposition \ref{dglie}, we saw that
$C^{*}(\A_{2},\A_{1})$ has a dg-Lie algebra structure.
We study a Maurer-Cartan equation in the dg-Lie algebra.
\begin{corollary}\label{maincorollary0}
The result of twisting
$\T_{H}=\A_{1}\oplus\A_{2}$ is also a twilled algebra
if and only if the curvature vanishes, or equivalently,
$H$ is a solution of a Maurer-Cartan equation,
$$
d_{\hat{\mu}_{2}}\what{H}+
\frac{1}{2}[\what{H},\what{H}]_{\hat{\mu}_{1}}=0.
\eqno{(MC)}
$$
The condition (MC) is equivalent with
\begin{equation}\label{defhami}
H(x)*_{1}H(y)+H(x)*_{2}y+x*_{2}H(y)=
H(H(x)*_{1}y+x*_{1}H(y))+H(x*_{2}y).
\end{equation}
\end{corollary}
\begin{proof}
We have $d_{\hat{\mu}_{2}}\what{H}=\hat{\mu}_{2}(\what{H}\ot 1)-
\what{H}\hat{\mu}_{2}+\hat{\mu}_{2}(1\ot\what{H})$ and
\begin{eqnarray*}
\frac{1}{2}[\what{H},\what{H}]_{\hat{\mu}_{1}}&=&
\frac{1}{2}\{\{\hat{\mu}_{1},\what{H}\},\what{H}\}\\
&=&\hat{\mu}_{1}(\what{H}\ot\what{H})-
\what{H}\hat{\mu}_{1}(1\ot\what{H})-
\what{H}\hat{\mu}_{1}(\what{H}\ot 1).
\end{eqnarray*}
This gives, for any $(a,x),(b,y)\in\T$,
\begin{multline*}
(d_{\hat{\mu}_{2}}\what{H}+
\frac{1}{2}[\what{H},\what{H}]_{\hat{\mu}_{1}})
((a,x),(b,y))=\\
H(x)*_{2}y-H(x*_{2}y)+x*_{2}H(y)+H(x)*_{1}H(y)
-H(H(x)*_{1}y+x*_{1}H(y)).
\end{multline*}
\end{proof}
\begin{definition}
Let $\A_{1}\Join\A_{2}$ be a twilled algebra and let
$H:\A_{2}\to\A_{1}$ a linear map.
We call the operator $H$ in (MC),
or equivalently, in (\ref{defhami})
a \textbf{Maurer-Cartan operator}.
A Maurer-Cartan operator is called \textbf{strong},
if it is a derivation with respect to
the multiplication $*_{2}$, i.e.,
$$
H(x*_{2}y)=x*_{2}H(y)+H(x)*_{2}y.
$$
\end{definition}
In Liu and coauthors \cite{LWX}, a Maurer-Cartan equation
in other dg-Lie algebra was studied.
The concept of strong solution is due to their work.
If $H$ is strong, then
the identity,
$H(x)*_{1}H(y)=H(H(x)*_{1}y+x*_{1}H(y))$,
automatically holds.
The strong Maurer-Cartan condition is equivalent with
$$
d_{\hat{\mu}_{2}}\what{H}=
\frac{1}{2}[\what{H},\what{H}]_{\hat{\mu}_{1}}=0.
$$
We easily obtain
\begin{corollary}\label{modelcoro}
If $H$ is a Maurer-Cartan operator, then
$$
x\times_{H}y:=H(x)*_{1}y+x*_{1}H(y)+x*_{2}y
$$
is an associative multiplication on $\A_{2}$.
\end{corollary}
\begin{proof}
When $H$ satisfies (MC), we have $\hat{\phi}^{H}_{2}=0$.
By Lemma \ref{5decomp}, we obtain
$\{\hat{\mu}^{H}_{2},\hat{\mu}^{H}_{2}\}=0$
which gives the associativity of $\hat{\mu}^{H}_{2}$.
The multiplication has the following form on $\A_{2}$,
$$
\hat{\mu}^{H}_{2}(x,y)=H(x)*_{1}y+x*_{1}H(y)+x*_{2}y.
$$
\end{proof}
We recall Rota-Baxter operators in Introduction.
\begin{example}
(Rota-Baxter operators of weight $q$.)
Let $\A$ be an associative algebra.
We recall the twilled algebra in Example \ref{nonabext}.
The multiplication of $\A\Join_{q}\A$ is defined by
\begin{equation}\label{ctaaq}
(a,x)*_{q}(b,y):=(ab,ay+xb+qxy),
\end{equation}
where $q\in\mathbb{K}$ (weight).
From (\ref{defhami}),
the Maurer-Cartan operators on $\A\Join_{q}\A$
satisfy the Rota-Baxter identity of weight $q$,
$$
R(x)R(y)=R(R(x)y+xR(y))+qR(xy).
$$
where we put $R:=H$.
Thus Rota-Baxter operators can be seen as examples
of Maurer-Cartan operators.\\
\indent
As an example of Rota-Baxter operator, we know
$$
R(f)(x):=f(qx)+f(q^{2}x)+f(q^{3}x)+...
\eqno{(\text{convergent})}
$$
where $R$ is defined on a certain algebra
of functions (see \cite{Rot2}).
\end{example}

\subsubsection{The cases of $\hat{\mu}_{2}=0$.}

Consider the cases of $\hat{\mu}_{2}=0$.
In this case, since $d_{\hat{\mu}_{2}}=0$,
the Maurer-Cartan equation simply has the form,
$[\hat{H},\hat{H}]_{\hat{\mu}_{1}}/2=0$,
or equivalently, (\ref{defhami}) reduces to the identity,
$$
H(x)*_{1}H(y)=H(H(x)*_{1}y+x*_{1}H(y)).
$$
Further, if $\A_{2}=\A_{1}$ as a canonical bimodule, then
$H$ is considered as a Rota-Baxter operator with weight zero.
\begin{definition}
(\cite{U1})
Let $\A$ be an associative algebra and let $M$ be an $\A$-bimodule.
A linear map $\pi:M\to\A$ is called a generalized Rota-Baxter
operator (of weight zero), if $\pi$ is a solution of
the identity,
\begin{equation}\label{(AP)}
\pi(m)\pi(n)=\pi(\pi(m)\cdot n+m\cdot\pi(n)),
\end{equation}
or equivalently, $[\hat{\pi},\hat{\pi}]_{\hat{\mu}}/2=0$,
where $m,n\in M$ and $\hat{\mu}$
is the associative structure of $\A\ltimes M$.
\end{definition}
A generalized Rota-Baxter
operator is obviously a (strong-)Maurer-Cartan operator.
Given a generalized Rota-Baxter operator $\pi:M\to\A$,
we have a twilled algebra $\A\Join M_{\pi}$ by the twisting
of $\A\ltimes M$ by $\pi$,
where $M_{\pi}$ is an associative subalgebra given by
Corollary \ref{modelcoro}.
The associative structure of $\A\Join M_{\pi}$
is the sum of two structures, $\hat{\mu}+\{\hat{\mu},\hat{\pi}\}$.
\begin{corollary}\label{satworo1}
Under the assumptions above, if $\pi_{1}$ is a second
generalized Rota-Baxter operator on $\A\ltimes M$, i.e.,
$[\hat{\pi}_{1},\hat{\pi}_{1}]_{\hat{\mu}}=0$,
then $H:=\pi_{1}-\pi$ is a Maurer-Cartan operator
on $\A\Join M_{\pi}$.
If $H$ is strong, then $\pi+tH$ is a one parameter family of
generalized Rota-Baxter operators for any $t\in\mathbb{K}$.
\end{corollary}
\begin{proof}
From assumptions, we have
$[\what{H},\what{H}]_{\hat{\mu}}/2=
-[\hat{\pi}_{1},\hat{\pi}]_{\hat{\mu}}$.
On the other hand,
since $d_{\hat{\mu}_{2}}(\cdot)=\{\{\hat{\mu},\hat{\pi}\},\cdot\}$,
we have
$$
d_{\hat{\mu}_{2}}\what{H}=\{\{\hat{\mu},\hat{\pi}\},\hat{\pi}_{1}\}=
[\hat{\pi},\hat{\pi}_{1}]_{\hat{\mu}}=
[\hat{\pi}_{1},\hat{\pi}]_{\hat{\mu}}.
$$
Simply, we obtain the condition (MC).
Thus Maurer-Cartan operators on $\A\Join M_{\pi}$ are given
as the difference of $\pi$ with generalized Rota-Baxter operators.
If $H$ is a strong Maurer-Cartan operator,
then $tH$ is also so for any $t\in\mathbb{K}$.
This implies the second part of the corollary.
\end{proof}

We recall in Section 3.1.2.
Let $\A$ be a finite dimensional associative algebra
and let $\A^{*}$ the dual space.
By a canonical adjoint action,
$\A$ acts on the dual space.
In this case, there are interesting similarities
in between generalized Rota-Baxter operators
and classical $r$-matrices.
We recall classical Yang-Baxter equation (CYBE).
There exists several equivalent definition of CYBE.
We recall the one of them.
CYBE is defined to be an operator identity
in the category of Lie algebras,
$$
\ [\tilde{r}(x),\tilde{r}(y)]=
\tilde{r}([\tilde{r}(x),y]+[x,\tilde{r}(y)])
$$
where $r$ is a two tensor in $\g\ot\g$
($\g$ is a finite dimensional Lie algebra),
$\tilde{r}:\g^{*}\to\g$ is the associated linear map,
$x,y$ are elements in the dual space $\g^{*}$
and the brackets in the right-hand side are adjoint actions.
The space of alternative tensors $\bigwedge^{*}\g$ has
a graded Lie algebra structure of Schouten bracket.
If $r$ is an element in $\g\wedge\g$,
then the Schouten bracket $[r,r]$
is in $\bigwedge^{3}\g$, and $[r,r]=0$ if and only if
$\tilde{r}$ satisfies CYBE above. Such a matrix $r$ is
called a triangular $r$-matrix. When $\g$ is a Lie algebroid,
a triangular $r$-matrix is a Poisson structure.
The notion of generalized Rota-Baxter operator can be seen as
an associative version of triangular $r$-matrices
and Poisson structures. We believe that
this picture is justified by the following example.
\begin{example}\label{rmatlab}
Let $\A$ be a 2-dimensional algebra generated by
$\left(
\begin{array}{cc}
0 & 1 \\
0 & 0
\end{array}
\right)$
and
$\left(
\begin{array}{cc}
1 & 0 \\
0 & 0
\end{array}
\right)$.
The dual space $\A^{*}$ is an $\A$-bimodule by adjoint action.
Thus we have a twilled algebra $\A\ltimes\A^{*}$.
Define a tensor $r$ by
$$
r:=
\left(
\begin{array}{cc}
0 & 1 \\
0 & 0
\end{array}
\right)
\wedge
\left(
\begin{array}{cc}
1 & 0 \\
0 & 0
\end{array}
\right).
$$
The tensor $r$ is identified with a map $\tilde{r}:\A^{*}\to\A$.
By direct computation, one can check that
the map is a generalized Rota-Baxter operator.\\
\indent
In general, if a 2-tensor $r\in\A\wedge\A$
satisfies Aguiar's multiplicative equation
(called an associative Yang-Baxter) in \cite{Ag,Ag3,Ag2},
$$
r_{13}r_{12}-r_{12}r_{23}+r_{23}r_{13}=0,
\eqno{(AYBE)}
$$
then $\tilde{r}:\A^{*}\to\A$ is a generalized Rota-Baxter
operator (see \cite{U1}).
Conversely, a skew symmetric generalized Rota-Baxter
operator satisfies (AYBE) above.
In non skewsymmetric cases,
there is a delicate difference between AYBE and
the generalized Rota-Baxter condition.\\
\indent
When $r$ is skewsymmetric,
the twisting by $r$ preserves the bilinear
pairing $(-|-)$ in Section 3.1.2.
Thus the associative structure $\hat{\mu}+\{\hat{\mu},\hat{r}\}$
satisfies the invariant condition in the sense of 3.1.2.
\end{example}

A Poisson structure is considered as a sheaf version
of triangular matrices. It is natural to ask what
is a sheaf version of Rota-Baxter operators.
We do not yet have an interesting solution.
We wish to find a Rota-Baxter operator on
the universal enveloping algebra of a Lie algebroid.
If there exists such a Rota-Baxter operator,
it is considered as an example of the sheaf version.

\subsection{The cases of $\phi_{1}\neq 0$ and $\phi_{2}=0$.}
In this case, $\T=\A_{1}\oplus\A_{2}$ is a quasi-twilled algebra.
However $\T_{H}=\A_{1}\oplus\A_{2}$ is not necessarily
a quasi-twilled algebra, because $\phi_{1}^{H}=\phi_{1}\neq 0$ and
$$
\hat{\phi}^{H}_{2}=d_{\hat{\mu}_{2}}\what{H}+
\frac{1}{2}[\what{H},\what{H}]_{\hat{\mu}_{1}}+
\frac{1}{6}\{\{\{\hat{\phi}_{1},\what{H}\},\what{H}\},\what{H}\}\neq 0.
$$
In general, the result of twisting have the forms,
\begin{eqnarray*}
\hat{\phi}^{H}_{1}&=&\hat{\phi}_{1},\\
\hat{\mu}^{H}_{1}&=&\hat{\mu}_{1}+\{\hat{\phi}_{1},\what{H}\},\\
\hat{\mu}^{H}_{2}&=&\hat{\mu}_{2}+
d_{\hat{\mu}_{1}}\what{H}+
\frac{1}{2}\{\{\hat{\phi}_{1},\what{H}\},\what{H}\},\\
\hat{\phi}^{H}_{2}&=&d_{\hat{\mu}_{2}}\what{H}+
\frac{1}{2}[\what{H},\what{H}]_{\hat{\mu}_{1}}+
\frac{1}{6}\{\{\{\hat{\phi}_{1},\what{H}\},\what{H}\},\what{H}\},
\end{eqnarray*}
Since $\hat{\mu}_{1}$ is not associative,
the derived bracket $[,]_{\hat{\mu}_{1}}$ does not satisfy
the graded Jacobi rule in general.
However the space $C^{*}(\A_{2},\A_{1})$ still
has a homotopy Lie algebra structure
$(d_{\hat{\mu}_{2}},[\cdot,\cdot]_{\hat{\mu}_{1}},
[\cdot,\cdot,\cdot]_{\hat{\phi}_{1}})$ in Section 3.2.
We consider a Maurer-Cartan equation
in this homotopy Lie algebra. The following
two corollaries are followed by the same manners
with Corollary \ref{maincorollary0} and Corollary \ref{modelcoro}.
\begin{corollary}
The result of twisting
$\T_{H}=\A_{1}\oplus\A_{2}$ is also a quasi-twilled algebra
if and only if it is a solution of twisted Maurer-Cartan equation,
$$
d_{\hat{\mu}_{2}}\what{H}+
\frac{1}{2}[\what{H},\what{H}]_{\hat{\mu}_{1}}+
\frac{1}{6}[\what{H},\what{H},\what{H}]_{\hat{\phi}_{1}}=0,
\eqno{(TMC)}
$$
or equivalently, for any $x,y\in\A_{2}$,
\begin{multline}\label{deftwisthami}
H(x)*_{1}H(y)+H(x)*_{2}y+x*_{2}H(y)=\\
H(H(x)*_{1}y+x*_{1}H(y))+H(x*_{2}y)+H(\phi_{1}(H(x),H(y))).
\end{multline}
\end{corollary}
\begin{corollary}\label{modelcoro2}
If $\T_{H}=\A_{1}\oplus\A_{2}$ is a quasi-twilled algebra, then
$$
x\times_{H,\phi_{1}}y:=\hat{\mu}^{H}_{2}(x,y)=
H(x)*_{1}y+x*_{1}H(y)+x*_{2}y+\phi_{1}(H(x),H(y)).
$$
is an associative multiplication on $\A_{2}$.
\end{corollary}

\begin{example}\label{twistbax}
(Twisted Rota-Baxter operators \cite{U1}.)
If $\hat{\mu}_{2}=0$, or $*_{2}$ is trivial,
then (\ref{deftwisthami}) is reduced to an identity:
$$
H(x)*_{1}H(y)=H(H(x)*_{1}y+x*_{1}H(y))+H(\phi_{1}(H(x),H(y))).
\eqno{(TRB1)}
$$
(TRB1) is equivalent with
$$
\frac{1}{2}[\what{H},\what{H}]_{\hat{\mu}_{1}}=
-\frac{1}{6}[\what{H},\what{H},\what{H}]_{\hat{\phi}_{1}}.
\eqno{(TRB2)}
$$
Such an operator $H$
is called a twisted Rota-Baxter operator (of weight zero).\\
\indent
As an example of twisted Rota-Baxter operators,
we know Reynolds operators in probability theory (\cite{Rot3}).
Let $\A$ be a certain functional algebra.
Define an operator $R:\A\to\A$ by
$$
R(f)(x):=\int^{\infty}_{0}e^{-t}f(x-t)dt.
$$
Then $R$ satisfies an identity,
$$
R(f)R(g)=R(R(f)g+fR(g))-R(R(f)R(g)).
$$
Such an operator is called a Reynolds operator.
The last term
$-R(R(f)R(g))=R\phi(R(f),R(g))$ can be seen as the cocycle term
of twisted Rota-Baxter identity.
Thus a Reynolds operator can be seen as a homotopy version of
Rota-Baxter operators of weight zero.
\end{example}
A Reynolds operator is used,
in the study of turbulent flow,
in order to induce a mean field model
of Navier-Stokes equation (so-called Reynolds equation).
One can easily verify that
if $R(f):=\overline{f}$ is the mean of $f$,
then the operator satisfies the identity above,
because an averaging operation satisfies the identities
$\overline{\overline{f}g}=\overline{f}\cdot\overline{g}
=\overline{f\overline{g}}$ and
$\overline{\overline{f}}=\overline{f}$ in general.
Unfortunately, we do not know an application
of our construction to Rota's theory.

\subsection{The cases of $\phi_{1}=0$ and $\phi_{2}\neq 0$}

In this case, $\hat{\phi}_{1}=\hat{\phi}^{H}_{1}=0$,
and thus $\hat{\mu}_{1}$ and $\hat{\mu}^{H}_{1}$
are both associative.
The twisted 4 substructures have the forms,
\begin{eqnarray*}
\hat{\mu}^{H}_{1}&=&\hat{\mu}_{1},\\
\hat{\mu}^{H}_{2}&=&\hat{\mu}_{2}+
d_{\hat{\mu}_{1}}\what{H},\\
\hat{\phi}^{H}_{2}&=&\hat{\phi}_{2}+d_{\hat{\mu}_{2}}\what{H}+
\frac{1}{2}[\what{H},\what{H}]_{\hat{\mu}_{1}}.
\end{eqnarray*}
Similar with Corollary \ref{maincorollary0}
and Corollary \ref{modelcoro},
we obtain the two corollaries below.
\begin{corollary}
The result of twisting
$\T_{H}=\A_{1}\oplus\A_{2}$ is a usual twilled algebra,
i.e., $\hat{\phi}^{H}_{2}=0$ if and only if
$H$ is a solution of the quasi-Maurer-Cartan equation,
$$
d_{\hat{\mu}_{2}}\what{H}+
\frac{1}{2}[\what{H},\what{H}]_{\hat{\mu}_{1}}=-\hat{\phi}_{2},
\eqno{(QMC)}
$$
or equivalently,
\begin{multline}\label{quasihami2}
H(x)*_{2}y+x*_{2}H(y)+H(x)*_{1}H(y)+\phi_{2}(x,y)=\\
H(H(x)*_{1}y+x*_{1}H(y))+H(x*_{2}y).
\end{multline}
\end{corollary}
\begin{corollary}
If $H$ satisfied (QMC), then
$\hat{\mu}^{H}_{2}$ is an associative structure
and defines an associative multiplication on $\A_{2}$ by
\begin{equation}\label{hdend}
x\times_{H,\phi_{2}}y:=
\hat{\mu}^{H}_{2}(x,y)=H(x)*_{1}y+x*_{1}H(y)+x*_{2}y.
\end{equation}
\end{corollary}
We consider a case of $\hat{\mu}_{2}=0$.
Then (QMC) and (\ref{quasihami2}) reduce to the identities,
respectively,
$$
\frac{1}{2}[\what{H},\what{H}]_{\hat{\mu}_{1}}=-\hat{\phi}_{2},
$$
and
\begin{eqnarray}\label{quasihami3}
H(x)*_{1}H(y)-H(H(x)*_{1}y+x*_{1}H(y))=-\phi_{2}(x,y).
\end{eqnarray}

Recall the quasi-twilled algebra $\A\oplus_{Q}\A$
in Example \ref{nonabext2}.
\medskip\\
\noindent
Claim.
Define a linear map $(a,x)\mapsto(\frac{q}{2}x,0)$ on $\A\oplus\A$.
Then its integral $e^{\widehat{q/2}}$ is an algebra isomorphism,
$$
e^{\widehat{q/2}}:
\A\Join_{q}\A\to\A\oplus_{Q}\A, \ \ Q=\frac{q^{2}}{4}.
$$
\begin{proof}
\begin{eqnarray*}
e^{\widehat{q/2}}((a,x)*_{q}(b,y))&=&
(ab+\frac{q}{2}ay+\frac{q}{2}xb+\frac{q^{2}}{2}xy,ay+xb+qxy)\\
&=&((a+\frac{q}{2}x)(b+\frac{q}{2}y)+\frac{q^{2}}{4}xy,ay+xb+qxy)\\
&=&(a+\frac{q}{2}x,x)*_{Q}(b+\frac{q}{2}y,y), \ \ Q=\frac{q^{2}}{4}.
\end{eqnarray*}
\end{proof}
If $Q=0$, then $\A\oplus_{Q=0}\A$ is the semi-direct product algebra.
Thus $\A\Join_{q}\A$ is isomorphic
with $\A\ltimes\A$ modulus $q^{2}$.\\
\indent
Now, the claim says that
$\A\Join_{q}\A$ is the result of twisting
of $\A\oplus_{Q}\A$ by $q/2$.
Let $(R(\A),\A)$ be the graph of $R$.
One can easily verify that
if $R$ is a $q$-Rota-Baxter operator,
then $\A\Join_{q}\A=\A\Join(R(\A),\A)$ is a second
twilled algebra decomposition.
By the twisting, we have a twilled algebra,
$\A\Join(R(\A)+\frac{q}{2}\A,\A)$,
$$
\A\Join(R(\A),\A)=\A\Join_{q}\A\overset{e^{\widehat{q/2}}}{\to}
\A\oplus_{q^{2}/4}\A=\A\Join(R(\A)+\frac{q}{2}\A,\A).
$$
\begin{example}
(Rota-Baxter operator mod $q^{2}$ \cite{E}).
Let $(\A,R)$ be a Rota-Baxter algebra.
We define a linear map $B:\A\to\A$ by
$B(\A):=R(\A)+\frac{q}{2}\A$.
Then the graph of $B$, $(B(\A),\A)$, is a subalgebra of
the quasi-twilled algebra $\A\oplus_{q^{2}/4}\A$.
This implies that $B$ is a solution of
$$
B(x)B(y)-B(B(x)y+xB(y))=-\frac{q^{2}}{4}xy.
$$
The right-hand term $q^{2}/4xy:=\phi_{2}(x,y)$
can be seen as the cocycle-term in (\ref{quasihami3}).
\end{example}

\section{Application.}

In this section, we will give a construction of
associative Nijenhuis operator.
First we recall basic properties of Nijenhuis operator.
A linear operator, $N:\A\to\A$,
is called an associative Nijenhuis operator,
if $N$ is a solution of
$$
N(x)N(y)=N(N(x)y+xN(y))-N^{2}(xy).
$$
In general, given a Nijenhuis operator,
$x\times_{N}y:=N(x)y+xN(y)-N(xy)$ is
a second associative multiplication and it is compatible
with the original multiplication.
Namely, $xy+tx\times_{N}y$
is a one parameter family of associative
multiplications for any $t\in\mathbb{K}$ (\cite{CGM}).
\medskip\\
\noindent
In the following,
we assume that $\A$ is an associative algebra,
$M$ is an $\A$-bimodule and
we denote the multiplication of $\A$ by $*_{\A}$.
\medskip\\
\indent
Let $\pi:M\to\A$ be a generalized Rota-Baxter operator, i.e.,
$\pi$ satisfies the identity,
\begin{equation}\label{generb}
\pi(m)*_{\A}\pi(n)=\pi(\pi(m)\cdot n+m\cdot \pi(n)).
\end{equation}
where $\cdot$ is the bimodule action of $\A$ on $M$ and $m,n\in M$.
We recall the twilled algebra $\A\Join M_{\pi}$ in Section 5.1.3.
The associative multiplication of $\A\Join M_{\pi}$ has the form
$$
(a,m)*(b,n)=
(a*_{\A}b+a\cdot_{\pi}n+m\cdot_{\pi}b,a\cdot n+m\cdot b+m\times_{\pi}n),
$$
where $\cdot_{\pi}$ means the bimodule action of $M_{\pi}$ on $\A$,
explicitly,
\begin{eqnarray*}
m\cdot_{\pi}b&:=&\pi(m)*_{\A}b-\pi(m\cdot b),\\
a\cdot_{\pi}n&:=&a*_{\A}\pi(n)-\pi(a\cdot n),
\end{eqnarray*}
and $m\times_{\pi}n$ is the associative multiplication of $M_{\pi}$,
explicitly,
$$
m\times_{\pi}n:=\pi(m)\cdot n+m\cdot \pi(n).
$$
Simply, we have $\pi(m\times_{\pi}n)=\pi(m)*_{\A}\pi(n)$.
\medskip\\
\indent
We consider a linear map $\Omega:\A\to M_{\pi}$.
The map $\Omega$ is a strong Maurer-Cartan operator
on a twilled algebra $M_{\pi}\Join\A$
if and only if
\begin{eqnarray}
\label{strong51}\Omega(a*_{\A}b)&=&a\cdot\Omega(b)+\Omega(a)\cdot b,\\
\label{strong52}\Omega(a)\times_{\pi}\Omega(b)&=&
\Omega(\Omega(a)\cdot_{\pi}b+a\cdot_{\pi}\Omega(b)),
\end{eqnarray}
or equivalently, $\Omega$ is a solution of
$$
d_{\hat{\mu}}\hat{\Omega}=
\frac{1}{2}[\hat{\Omega},\hat{\Omega}]_{\{\hat{\mu},\hat{\pi}\}}=0.
$$
We give the main result of this section.
\begin{proposition}\label{constN}
Let $\Omega:\A\to M_{\pi}$ be a strong Maurer-Cartan operator.
\begin{description}
	\item[1.] Then a composition map $N:=\pi\Omega$ is an associative
	Nijenhuis operator on $\A$. Namely $N$ satisfies the condition
	$$
	N(a)*_{\A}N(b)=N(N(a)*_{\A}b+a*_{\A}N(b))-NN(a*_{\A}b)
	$$
	for any $a,b\in\A$.
\end{description}
The pair of $(\pi,N)$ is compatible in the following sense.
\begin{description}
	\item[2.] The composition $N\pi:M\to\A$ is a second
	generalized Rota-Baxter operator.
	\item[3.] The operators $\pi$ and $N\pi$ are compatible, i.e.,
	$$
	[\hat{\pi},\widehat{N\pi}]_{\hat{\mu}}=0.
	$$
	This implies that $N\pi$ is strong as a Maurer-Cartan operator
	and $\pi+tN\pi$ $t\in\mathbb{K}$ is a one parameter family
	of generalized Rota-Baxter operators.
\end{description}
\end{proposition}
\begin{proof}
1. Applying $\pi$ to (\ref{strong52}), we have
$$
\pi\Omega(a)*_{\A}\pi\Omega(b)=
\pi\Omega(\Omega(a)\cdot_{\pi}b+a\cdot_{\pi}\Omega(b)).
$$
In the right-hand side,
$$
\Omega(a)\cdot_{\pi}b+a\cdot_{\pi}\Omega(b)=
\pi\Omega(a)*_{\A}b-\pi(\Omega(a)\cdot b)+
a*_{\A}\pi\Omega(b)-\pi(a\cdot\Omega(b)).
$$
From (\ref{strong51}), we have
$$
\Omega(a)\cdot_{\pi}b+a\cdot_{\pi}\Omega(b)=
\pi\Omega(a)*_{\A}b+a*_{\A}\pi\Omega(b)-\pi\Omega(a*_{\A}b).
$$
Thus we obtain the desired condition,
$$
\pi\Omega(a)*_{\A}\pi\Omega(b)=
\pi\Omega(\pi\Omega(a)*_{\A}b+a*_{\A}\pi\Omega(b))-
\pi\Omega\pi\Omega(a*_{\A}b).
$$
2. We put $a:=\pi(m)$ and $b:=\pi(n)$
for any $m,n\in M$.
Then, by the Nijenhuis condition of $\pi\Omega$, we have
\begin{equation}\label{2mainp}
\pi\Omega\pi(m)*_{\A}\pi\Omega\pi(n)=
\pi\Omega(\pi\Omega\pi(m)*_{\A}\pi(n)+\pi(m)*_{\A}\pi\Omega\pi(n))
-\pi\Omega\pi\Omega(\pi(m)*_{\A}\pi(n)).
\end{equation}
From the identity (\ref{generb}), we have
\begin{eqnarray*}
\pi\Omega\pi(m)*_{\A}\pi(n)&=&
\pi(\pi\Omega\pi(m)\cdot{n}+\Omega\pi(m)\cdot\pi(n)),\\
\pi(m)*_{\A}\pi\Omega\pi(n)&=&
\pi(\pi(m)\cdot\Omega\pi(n)+m\cdot\pi\Omega\pi(n)),
\end{eqnarray*}
and from the derivation rule, we have
$$
\pi\Omega\pi\Omega(\pi(m)*_{\A}\pi(n))=
\pi\Omega\pi(\Omega\pi(m)\cdot\pi(n)+\pi(m)\cdot\Omega\pi(n)).
$$
Thus (\ref{2mainp}) has the form,
\begin{multline*}
\pi\Omega\pi(m)*_{\A}\pi\Omega\pi(n)=
\pi\Omega\pi(\pi\Omega\pi(m)\cdot{n}+\Omega\pi(m)\cdot\pi(n)+
\pi(m)\cdot\Omega\pi(n)+m\cdot\pi\Omega\pi(n))-\\
\pi\Omega\pi(\Omega\pi(m)\cdot\pi(n)+\pi(m)\cdot\Omega\pi(n))=\\
\pi\Omega\pi(\pi\Omega\pi(m)\cdot{n}+m\cdot\pi\Omega\pi(n)),
\end{multline*}
this is the desired result.\\
3. It is obvious that
$\widehat{\pi\Omega\pi}=\hat{\pi}\hat{\Omega}\hat{\pi}$.
We have $[\hat{\pi},\widehat{\pi\Omega\pi}]_{\hat{\mu}}=$
\begin{multline}\label{pipiopi}
\{\{\hat{\mu},\hat{\pi}\},\hat{\pi}\hat{\Omega}\hat{\pi}\}=
\{\hat{\mu}(\pi\ot 1)+\hat{\mu}(1\ot\hat{\pi})-\hat{\pi}\hat{\mu},
\hat{\pi}\hat{\Omega}\hat{\pi}\}=\\
\hat{\mu}(\hat{\pi}\ot\hat{\pi}\hat{\Omega}\hat{\pi})-
\hat{\pi}\hat{\Omega}\hat{\pi}\hat{\mu}(\hat{\pi}\ot 1)+
\hat{\mu}(\hat{\pi}\hat{\Omega}\pi\ot\hat{\pi})-
\hat{\pi}\hat{\Omega}\hat{\pi}\hat{\mu}(1\ot\hat{\pi})\\
-\hat{\pi}\hat{\mu}(\hat{\pi}\hat{\Omega}\hat{\pi}\ot 1)
-\hat{\pi}\hat{\mu}(1\ot\hat{\pi}\hat{\Omega}\hat{\pi}),
\end{multline}
where $\hat{\pi}\hat{\pi}=0$ is used.
From the generalized Rota-Baxter condition,
$[\hat{\pi},\hat{\pi}]_{\hat{\mu}}/2=
\hat{\mu}(\hat{\pi}\ot\hat{\pi})
-\hat{\pi}\hat{\mu}(\hat{\pi}\ot{1})-\hat{\pi}\hat{\mu}(1\ot\hat{\pi})=0$,
we have
\begin{multline}\label{pipiopi2}
(\ref{pipiopi})=
\hat{\mu}(\hat{\pi}\ot\hat{\pi}\hat{\Omega}\hat{\pi})-
\hat{\pi}\hat{\Omega}\hat{\mu}(\hat{\pi}\ot\hat{\pi})+
\hat{\mu}(\hat{\pi}\hat{\Omega}\pi\ot\hat{\pi})
-\hat{\pi}\hat{\mu}(\hat{\pi}\hat{\Omega}\hat{\pi}\ot 1)
-\hat{\pi}\hat{\mu}(1\ot\hat{\pi}\hat{\Omega}\hat{\pi})=\\
-\hat{\pi}\hat{\Omega}\hat{\mu}(\hat{\pi}\ot\hat{\pi})+
\hat{\mu}(\hat{\pi}\hat{\Omega}\pi\ot\hat{\pi})
-\hat{\pi}\hat{\mu}(\hat{\pi}\hat{\Omega}\hat{\pi}\ot 1)
+\hat{\pi}\hat{\mu}(\hat{\pi}\ot\hat{\Omega}\hat{\pi})=\\
-\hat{\pi}\hat{\Omega}\hat{\mu}(\hat{\pi}\ot\hat{\pi})
+\hat{\pi}\hat{\mu}(\hat{\Omega}\hat{\pi}\ot\hat{\pi})
+\hat{\pi}\hat{\mu}(\hat{\pi}\ot\hat{\Omega}\hat{\pi}).
\end{multline}
Since $\hat{\Omega}$ is a derivation with respect to $\hat{\mu}$,
the last equation of (\ref{pipiopi2}) is zero.
\end{proof}
\begin{example}\label{exN}
We put $\A:=C^{1}([0,1])$ and $M:=C^{0}([0,1])$.
We assume a canonical bimodule action of $\A$ on $M$.
An integral operator is a Rota-Baxter operator with weight zero.
$$
\pi:M\to\A, \ \ \pi(f)(x):=\int^{x}_{0}dtf(t).
$$
Then a derivation from $\A$ to $M_{\pi}$,
$$
\Omega(f)(x):=\omega(x)\frac{df}{dx}(x)=\omega(x)f^{\p}(x), \ \
\omega(x)\in C^{0}([0,1])
$$
is a strong Maurer-Cartan operator.
The induced Nijenhuis operator on $\A$ is
$$
N(f)(x)=\int^{x}_{0}\omega(t)f^{\p}(t)dt.
$$
\end{example}
\begin{proof}
We only check the condition (\ref{strong52}).
For any $f,g\in\A$,
$$
\Omega(f)\cdot_{\pi}g=\pi\Omega(f)g-\pi(\Omega(f)g)
=\int^{x}_{0}dt\omega(t)f^{\p}(t)g(x)-
\int^{x}_{0}dt\omega(t)f^{\p}(t)g(t).
$$
We have
\begin{eqnarray*}
\Omega(\Omega(f)\cdot_{\pi}g)&=&
\int^{x}_{0}dt\omega(t)f^{\p}(t)\omega(x)g^{\p}(x),\\
\Omega(f\cdot_{\pi}\Omega(g))&=&
\omega(x)f^{\p}(x)\int^{x}_{0}dt\omega(t)g^{\p}(t).
\end{eqnarray*}
On the other hand,
\begin{eqnarray*}
\Omega(f)\times_{\pi}\Omega(g)&=&\omega(x)f^{\p}(x)\times_{\pi}\omega(x)g^{\p}(x)\\
&=&\int^{x}_{0}dt\omega(t)f^{\p}(t)\omega(x)g^{\p}(x)+
\omega(x)f^{\p}(x)\int^{x}_{0}dt\omega(t)g^{\p}(t).
\end{eqnarray*}
Thus we obtain the desired condition.
\end{proof}
We consider two examples in noncommutative cases.
In the proof of Example \ref{exN},
we used the commutativity of only $\omega$.
Hence if $\omega$ is $1$ or a central element,
then the similar proof holds over noncommutative setting.
\begin{example}
Let $\A$ be an associative algebra and let $\A[[\nu]]$
an algebra of formal series with coefficients in $\A$.
The multiplication on $\A[[\nu]]$ is defined by
$$
a_{i}\nu^{i}*b_{j}\nu^{j}:=a_{i}b_{j}\nu^{i+j},
\ \ a_{i},b_{j}\in\A,
$$
where $\sum$ is omitted.
We define a formal integral operator,
$$
\int d\nu a_{i}\nu^{i}:=\frac{1}{i+1}a_{i}\nu^{i+1}, \ \ a_{i}\in\A.
$$
The integral operator is a Rota-Baxter operator with weight zero.
The formal derivation operator is a strong Maurer-Cartan operator
$$
\Omega(a_{i}\nu^{i}):=z_{k}\nu^{k}\frac{d}{d\nu}(a_{i}\nu^{i})
:=iz_{k}a_{i}\nu^{i+k-1}, \ \ z_{k}\in Z(\A).
$$
Here $Z(\A)$ is the space of central elements.
The induced Nijenhuis operator is
$$
N(a_{i}\nu^{i}):=\frac{i}{i+k}z_{k}a_{i}\nu^{i+k}.
$$
\end{example}

\begin{example}\label{ewn}
Let $W\langle x,\partial_{x}\rangle$ be the Weyl algebra.
Define a formal integral operator by,
for the normal basis of the Weyl algebra,
$$
\int{dx\partial_{x}^{i}*x^{j}}:=\frac{1}{1+j}\partial_{x}^{i}*x^{j+1},
\ \ i,j\ge 0.
$$
Then the integral operator is a Rota-Baxter operator
with weight zero (see \cite{U1}).
We put $\Omega:=i_{\partial_{x}}$.
Then $\Omega$ is a strong Maurer-Cartan operator.
Thus the composition map
$$
N(u):=\int{dx}\Omega(u)=\int{dx}[\partial_{x},u]
$$
is a Nijenhuis operator on $W\langle x,\partial_{x}\rangle$.
Since an arbitrary element $u$ has the form of
$u:=k_{ij}\partial_{x}^{i}*x^{j(j\neq 0)}+k_{i}\partial_{x}^{i}+k$,
we have $N(u)=k_{ij}\partial_{x}^{i}*x^{j (j\neq 0)}$.
Thus $N$ is a projection onto the space of elements
of the form $k_{ij}\partial_{x}^{i}*x^{j(j\neq 0)}$.
\end{example}

\end{document}